\documentclass[a4paper, 11pt]{article}
\usepackage{geometry}
\usepackage[colorlinks]{hyperref}
\usepackage{comment}
\usepackage{amsmath, amsfonts, amsthm}
\usepackage{cite}
\usepackage{lmodern}
\usepackage{tikz}
\usetikzlibrary{patterns}
\usepackage{sgamevar}
\usepackage{enumitem}  \setenumerate[0]{label=(\roman*)}
\usepackage{yfonts}

\allowdisplaybreaks

\DeclareMathOperator*{\isom}{Isom}
\DeclareMathOperator*{\bij}{Bij}
\DeclareMathOperator*{\Aut}{Aut}
\newcommand{\id}{\text{id}}
\DeclareMathOperator*{\argmax}{arg\,max}
\DeclareMathOperator*{\argmin}{arg\,min}

\newcommand{\set}[1]{\left\{#1\right\}}
\newcommand{\reals}{\mathbb{R}}
\newcommand{\posreals}{\mathbb{R}^+}

\let\oldproofname=\proofname
\renewcommand{\proofname}{\rm\bf{\oldproofname}}

\usepackage{setspace, tocloft}

\newtheoremstyle{paper}
{3pt}
{3pt}
{}
{}
{\bfseries}
{:}
{.5em}
{}

\theoremstyle{paper}
\newtheorem{theorem}{Theorem}[section]
\newtheorem{proposition}[theorem]{Proposition}
\newtheorem{lemma}[theorem]{Lemma}
\newtheorem{corollary}[theorem]{Corollary}

\newtheorem{definition}[theorem]{Definition}
\newtheorem{example}[theorem]{Example}

\begin{document}
    \title{Notions of Anonymity, Fairness and Symmetry \\ for Finite Strategic-Form Games}
    \author{Nicholas Ham\footnote{{\it Email:} \href{mailto:contact@n-ham.com}{\tt contact@n-ham.com}}}
	\maketitle

    \begin{abstract}
    		In this paper we survey various notions of anonymity and symmetry for finite strategic-form games present in relevant literature, and discuss notions of fairness; show that game bijections and game isomorphisms form groupoids; introduce matchings as a convenient characterisation of strategy triviality; and outline how to construct and partially order parameterised (symmetric) games with examples that range all combinations of surveyed symmetry notions, which when combined with other results in this paper gives the precise relationship between the various symmetry notions.
    		
    		
    \end{abstract}
    
    \begin{center}
    	\tableofcontents
    \end{center}
    
    \section{Introduction} \label{sec:intro}
The notion of a game being fair may be made more precise with the concept of symmetry. Broadly speaking we will consider a game fair when the players are indifferent between which position they play, however there are several distinct notions of symmetry that are possible which lead to variations in structure and fairness. For example, the players may or may not care about the arrangement of their opponents. 

Note that this paper does not survey the literature on notions of symmetry, though the reader may find it a useful reference if undertaking such an endeavour. However this paper does attempt to survey the numerous notions of symmetry for finite strategic-form games that are present in the literature, whilst also filling various holes and opening several further directions of research in the area. This is important to our understanding of the theory of symmetric games and fairness, which is fundamental when it comes to the theory of games, artificial intelligence, biology, computer science, economic theory, evolution, legal systems, logic, philosophy, political science, along with social choice and voting theory to name just a few examples. See for example Arrow's impossibility theorem \cite{arrow1950difficulty, arrow2012social}, biological warfare \cite{williams1989unit, gold2011unit}, financial contagion \cite{dungey2005contagion}, forced breakups and relationships \cite{chantler2009forced, ouattara1998forced}, forced psychiatric treatments \cite{minkowitz2006united}, gaslighting \cite{sweet2019sociology, ahern2018institutional}, human misery \cite{margolis2003misery, curato2019democracy}, human trafficking \cite{aronowitz2009human, shelley2010human, hesse2011young}, incentive theory \cite{laffont1993theory, laffont2009theory}, psychological warfare \cite{farago1941german, linebarger2015psychological} and theory of capitalism \cite{harriss2006poverty} and democracy \cite{delli2002internet, walker1966critique}.

Symmetry and fairness in the context of games was first explored by von Neumann and Morgenstern \cite{VNM}, outlining what we will later refer to as our label-dependent framework in which player permutations act on strategy profiles, consequently requiring all players have the same strategy labels. Soon after Nash \cite{NashNCG} famously showed that symmetric games have at least one symmetric mixed strategy Nash equilibrium, while more recently Hofbauer et al. \cite{hofbauer2002differential} and Cheng et al. \cite{CRVWSym} showed that fully symmetric $2$-strategy games have at least one pure strategy Nash equilibrium. More on equilibria for symmetric games can be found in \cite{hefti2017equilibria}. Notions of symmetry and equivalence also appear in \cite{HarsanyiSelten}.

Notably the term fair has not really been used in the context of non-zero-sum strategic-form games, and quite a lot of people seem unwilling to discuss or consider other people's viewpoints on the topic of fairness, often doing their best to gain enough power to dictate their views on other people, too often disguising the intentions of their rules. However the term fair did appear as early as the 1940s in the context of zero-sum games, which are a subclass of strategic-form games, including von Neumann and Morgenstern \cite[17.11, 28.1, 28.2]{VNM} and for example \cite{gelbaum1959symmetric, tucker1962comb}. There will be a discussion in Subsection \ref{subsec:fairnessdiscussion} defending the author's use of the term fair in the context of symmetric strategic-form games, though note it is an incredibly complicated and intricate topic, with few to no objectively unambiguously correct answers, but may hopefully help contribute towards giving people an actual choice in life without nasty ways of creating victims, turning victims against each other and having people submit to slavery-like conditions through poverty etc..

Under the theme of anonymity rather than symmetry, Brandt et al. \cite{brandt2009symmetries} and Tohm\'e et al. \cite{tohme2019structural} examined label-dependent notions such as where players are indifferent between who plays which strategy, and where players do not distinguish between their opponents.

A number of people have examined notions of symmetry which may not be captured inside our label-dependent framework, see for example Nash \cite{NashNCG}, Shapley \cite{shapley1960symmetric}, Peleg et al. \cite{peleg1999canonical}, Sudh\"{o}lter et al. \cite{sudholter2000canonical} and Stein \cite{NoahXE}. In order to discuss and analyse such notions we will need to make a detour to examine morphisms between games, the complexity of which has been investigated by Gabarr\'{o} et al. \cite{IsoComplexity}. Inside what will later be referred to as our label-independent framework game automorphisms act on strategy profiles, which also allows players to have distinct strategy labels.

We begin in Section \ref{sec:background} by reviewing several mathematical concepts that will play an important role throughout our analysis. In Section \ref{sec:labdepnotions} we survey various label-dependent notions of anonymity, symmetry and fairness. 

In Section \ref{sec:morphisms} we review game morphisms while showing that game bijections and game isomorphisms form groupoids, which appears to be missing from relevant literature, and introduce matchings as a convenient characterisation of strategy triviality. 

Finally, in Section \ref{sec:labindnotions} we survey various label-independent notions of symmetry, discuss how to classify a given game, and outline how to construct and partially order parameterised symmetric games with examples that range over various classes.

    \section{Background} \label{sec:background}
	Let $N = \{1, \ldots, n\}$ where $n \geq 2$ and let $\{A_i: i \in N\}$ be a collection of non-empty sets. To simplify notation:
	\begin{enumerate}
		\item We denote the Cartesian product of $\{A_i: i \in N\}$, ie. $\times_{i \in N}A_i$, as $A$;
		\item For each $i \in N$ we denote $\times_{j \in N-\{i\}} A_j$ as $A_{-i}$;
		\item For each $s \in \times_{i \in N}A_i$ and $i \in N$ we denote the element of $A_i$ used in $s$, which is position $i$ of $s$, as $s_i$;
		\item For each $i \in N$ and $s = (s_1, \ldots, s_{i-1}, s_i, s_{i+1}, \ldots, s_n) \in A$ we denote $(s_1, \ldots, s_{i-1}, s_{i+1}, \ldots, s_n) \in A_{-i}$ as $s_{-i}$;
		\item For each $s_i \in A_i$ and $s_{-i} = (s_1, \ldots, s_{i-1}, s_{i+1}, \ldots, s_n) \in A_{-i}$ we denote $(s_1, \ldots, s_{i-1}, s_i, s_{i+1}, \ldots, s_n) \in \times_{i \in N} A_i$ as $(s_i, s_{-i})$ or $s$, whichever is contextually convenient;
		\item For each $k \in N$ and distinct $i_1, \ldots, i_k \in N$ we denote $N - \{i_1, \ldots, i_k\}$ as $N_{-i_1 - \ldots - i_k}$, for example $N_{-1 - \ldots - k} = \set{k+1, \ldots, n}$;
		\item For each $k \in N-\{n\}$ and distinct $i_1, \ldots, i_k \in N$ we denote
		\begin{align*}
			{\left({(A_{-i_1})}_{-i_2}\ldots\right)}_{-i_k} &= A_{{-i_1}_{{-i_2}_{\ldots_{-i_k}}}} \\
			&= \times_{i \in N-\{i_1, \ldots, i_k\}}A_i	
		\end{align*}
		as $A_{-i_1 - \ldots - i_k}$; and
		\item For each $s \in A$, $k \in N-\{n\}$ and distinct $i_1, \ldots, i_k \in N$ we denote \[{\left({(s_{-i_1})}_{-i_2}\ldots\right)}_{-i_k} = s_{{-i_1}_{{-i_2}_{\ldots_{-i_k}}}} \in A_{-i_1 - \ldots - i_k} = \times_{i \in N-\{i_1, \ldots, i_k\}}A_i\] as $s_{-i_1 - \ldots - i_k}$;
		\item For each $k \in N-\{n\}$, distinct $i_1, \ldots, i_k, j \in N$ and $s_{-i_1 - \ldots - i_k} \in A_{-i_1 - \ldots - i_k}$ we denote the element of $A_j$ used in $s_{-i_1 - \ldots - i_k}$ as ${(s_{-i_1 - \ldots - i_k})}_j$, note if $j>\min\set{i_l: l\in\{1, \ldots, k\}}$ then this will not be position $j$ of the $(n-k)$-tuple $s_{-i_1 - \ldots - i_k}$.
	\end{enumerate}
	
	A \textit{relation} on $\set{A_i: i \in N}$ is a subset $R$ of their Cartesian product $\times_{i \in N}A_i$. Let $i \in N$, we say that $R$ is \textit{$i$-total} when for each $s_i \in A_i$ there exists $s_{-i} \in A_{-i}$ such that $(s_i, s_{-i}) \in R$, and \textit{$i$-unique} when $(s_i, s_{-i}), (s_i, s_{-i}') \in R$ implies $s_{-i} = s_{-i}'$.
	
	Given sets $X$ and $Y$, a \textit{function from $X$ to $Y$} is a functional left-total binary relation $f \subseteq X\times{Y}$. We denote the set of functions from $X$ to $Y$ as $Y^X$. Let $f \in Y^X$:
	\begin{enumerate}
		\item Since $f$ is functional, for each $x \in X$ we may denote by $f(x)$ the unique element in $Y$ such that $(x, f(x)) \in f$;
		\item The \textit{image} of $f$ is the set $\set{f(x): x \in X}$, which we denote as $f(X)$; and
		\item If $Y=X$ then the function that maps each element of $X$ to itself acts as an identity under composition, it is typically referred to as \textit{the identity function} and denoted as $\id_X$. 
	\end{enumerate}
	
	A function $f \in Y^X$ is referred to as:
	\begin{enumerate}
		\item \textit{injective}, or as an \textit{injection}, when for each $x, x' \in X$, $f(x) = f(x')$ if and only if $x = x'$;
		\item \textit{surjective}, or as a \textit{surjection}, when $f(X) = Y$; and
		\item \textit{bijective}, or as a \textit{bijection}, when it is both injective and surjective. 
	\end{enumerate}
	
	Note that the bijections from a set to itself form a group under composition.
	
	A function $f \in \reals^{\reals}$ is referred to as:
	\begin{enumerate}
		\item a \textit{strictly increasing function} when for each $x, x' \in \mathbb{R}$, $x < x'$ if and only if $f(x) < f(x')$; and
		\item a \textit{positive linear transformation} when there exists $\alpha \in \posreals$ and $\beta \in \reals$ such that $f(x) = \alpha{x} + b$ for all $x \in X$.
	\end{enumerate}
	
	Note that the strictly increasing functions are a subgroup of the bijections from $\reals$ to $\reals$, and that the positive linear transformations are a subgroup of the strictly increasing functions. For proofs see \cite[Propositions 2.2.4 and 2.2.6]{ham2011honoursthesis}.
	
	A \textit{strategic-form game}, or just \textit{game} when contextually unambiguous, consists of a set $N = \{1, \ldots, n\}$ of $n \geq 2$ \textit{players}, or \textit{player names}, and for each player $i \in N$, a non-empty set $A_i$ of \textit{strategies} and a \textit{utility function} $u_i:A\rightarrow\mathbb{R}$, where $A$ denotes the set of \textit{strategy profiles} $\times_{i \in N}A_i$. We denote such a game as the triple $(N, A, u)$, where $u = (u_i)_{i \in N}$. If there exists $m \in \mathbb{Z}^+$ such that $|A_i| = m$ for all $i \in N$ then $(N, A, u)$ is called an $m$-\textit{strategy} game. A game $(N, A, u)$ is \textit{finite} when both $N$ is finite and $A_i$ is finite for all $i \in N$. 
	
	In this paper we will only concern ourselves with finite games, consequently all player sets and pure strategy sets are implicitly finite. It is noted in Mas-Collel, Whinston and Green \cite[Proposition 3.C.1]{Mas1995microeconomic} that for a set $X$ any rational preference relation may be described by some utility function, it is also worth checking Sections 3.B and 3.C, especially the definition of a rational preference relation (Definition 3.B.1). This would be a useful place to begin when trying to explore notions of symmetry and fairness for non-finite games (so a [possibly] non-finite number of players and/or [possibly] non-finite strategy sets, with at least one non-finite set involved). 
	
	Note that the strategy profiles, and consequently also the utility functions, of a game have an implicit ordering of the players. We refer to the place of each player in this order as their \textit{position}. For games when the player names are $\{1, \ldots, n\}$, unless otherwise specified, the names and positions coincide.
	
	A game may be displayed pictorially as a list of matrices. We list the strategies from players $n-1$ and $n$ along the rows and columns respectively (or from the players in positions $n-1$ and $n$ where the player names are not $\{1, \ldots, n\}$), and for games with more than two players have a separate matrix for each strategy combination of the remaining players $\{1, \ldots, n-2\}$. Each strategy profile $s \in A$ corresponds to a unique cell in one of the matrices where the payoffs are written in the form $\bigl(u_i(s)\bigr)_{i \in N}$. For an example, see Example \ref{fullsymeg}.
	
	\begin{example} \label{fullsymeg} 3-player 2-strategy game.
		\begin{center}
  		\begin{game}{2}{2}[$(a,,)$]
    			     \>  $a$      \>  $b$      \\
    			$a$  \>  $1,1,1$  \>  $2,2,3$  \\
    			$b$  \>  $2,3,2$  \>  $4,5,5$  
  		\end{game}
  		\hspace*{5mm}
  		\begin{game}{2}{2}[$(b,,)$]
       			 \>  $a$      \>  $b$      \\
    			$a$  \>  $3,2,2$  \>  $5,4,5$  \\
    			$b$  \>  $5,5,4$  \>  $6,6,6$ 
  		\end{game}
		\end{center}
		We find the payoff to player $3$ for the strategy profile $(b, b, a) \in A$ as follows: reading the strategy profile from left to right, player $1$ has chosen the second matrix, player $2$ has chosen the second row and player $3$ has chosen the first column, the third value of which is the payoff to player $3$. Hence $u_3(b,b,a) = 4$.
	\end{example}
	
	The reader should note that the \textit{usual convention} in most of the game theory literature is to have players $1$ and $2$ along the rows and columns respectively, and for games with more than two players have a separate matrix for each strategy combination of the remaining players $\{3, \ldots, n\}$. The author considers the usual convention objectively inferior to the convention used in this paper, primarily when finding the payoffs for a given strategy profile. 
		
	The normal and convenient convention is to read an $n$-tuple left-to-right, when reading a strategy profile left-to-right to find the payoffs for players:
	\begin{enumerate}
		\item using the usual convention first one finds the correct row and column in the first matrix, then one finds the correct matrix while trying to recall what the correct row and column are, which is incredibly tedious, frustrating and error-prone; whereas
		\item using the convention in this paper one first finds the correct matrix, then one finds the correct row and column, which also has one indexing the payoff matrices using the normal convention for matrices. 
	\end{enumerate} 
	
	A possible solution with the usual convention when finding the payoffs for a given strategy profile is to read strategy profiles left-to-right from player $3$ through to player $n$, and then players $1$ and $2$. The author still finds this less efficient and more tedious, frustrating and error-prone than simply changing to the convention used in this paper.
	
	Given a set $X$, we denote the set of probability distributions over $X$ as $\Delta(X)$, ie. $\Delta(X) = \{\sigma \in [0,1]^X: \Sigma_{x \in X}\sigma(x) = 1\}$. Need to clean this up.
	
	Given a game $\Gamma = (N, A, u)$, for each player $i \in N$, the \textit{mixed strategy set for player $i$} is the set of probability distributions over $A_i$, ie. $\Delta(A_i) = \{\sigma_i \in [0,1]^{A_i}: \Sigma_{s_i \in A_i} \sigma_i(s_i) = 1\}$. The set of \textit{mixed strategy profiles} is the Cartesian product of the players' mixed strategy sets, ie. $\times_{i \in N}\Delta(A_i)$, note that this is not the same as the probability distributions over $A$, ie. $\times_{i \in N}\Delta(A_i) \neq \Delta(A)$. To simplify notation we shall denote $\times_{i \in N}\Delta(A_i)$ as $\nabla(A)$. Given our notation from earlier, for each $i \in N$ we have $\nabla(A)_{-i} = \times_{j \in N-\{i\}}\Delta(A_j)$. 
	
	Define $f:\nabla(A)\rightarrow\Delta(A)$ where for each $\sigma = (\sigma_1, \ldots, \sigma_n) \in \nabla(A)$, we let $f(\sigma)(s) = f(\sigma_1, \ldots, \sigma_n)(s) = \prod_{i \in N}\sigma_i(s_i)$ for all $s \in A$. 
	
	Note that for each $k \in N-\{n\}$, distinct $i_1, \ldots, i_k \in N$, $\sigma \in \nabla(A)$ and $s \in A$:
	\begin{align*}
		f(\sigma_{-i_1 - \ldots - i_{k-1}})(s_{-i_1 - \ldots - i_{k-1}}) &= \prod_{j \in N_{-i_1 - \ldots - i_{k-1}}}{(\sigma_{-i_1 - \ldots - i_{k-1}})}_j({(s_{-i_1 - \ldots - i_{k-1}})}_j) \\
		&= \sigma_{i_k}(s_{i_k})\prod_{j \in N_{-i_1 - \ldots - i_k}}{(\sigma_{-i_1 - \ldots - i_k})}_j({(s_{-i_1 - \ldots - i_k})}_j) \\
		&= \sigma_{i_k}(s_{i_k})f(\sigma_{-i_1 - \ldots - i_k})(s_{-i_1 - \ldots - i_k}). \\
	\end{align*}
	
	\begin{proposition}
		The function $f$ defined above satisfies:
		\begin{enumerate}
			\item For each $\sigma = (\sigma_1, \ldots, \sigma_n) \in \nabla(A)$, $f(\sigma) = f(\sigma_1, \ldots, \sigma_n) \in \Delta(A)$ (ie. $f$ is well-defined);
			\item For each $\sigma, \sigma' \in \nabla(A)$, $f(\sigma) = f(\sigma')$ if and only if $\sigma = \sigma'$ (ie. $f$ is injective); and
			\item $f\left(\nabla(A)\right) \subset \Delta(A)$ (ie. $f$ is not surjective).
		\end{enumerate}
		
		\begin{proof}
			\begin{enumerate}
				\item First note that:
				\begin{align*}
					\sigma_i(s_i) &\geq 0 \text{ for all } i \in N, s_i \in A_i \\
					\Rightarrow f(\sigma)(s) = f(\sigma_1, \ldots, \sigma_n)(s) = \prod_{i \in N} \sigma_i(s_i) &\geq 0 \text{ for all } s \in A.
				\end{align*}
				Also we have:
				\begin{align*}
					\sum_{s \in A}f(\sigma)(s) = \sum_{s \in A}f(\sigma_1, \ldots, \sigma_n)(s) &= \sum_{s \in A}\prod_{i \in N}\sigma_i(s_i) \\
					&= \sum_{s \in A}\left[\sigma_1(s_1)f(\sigma_{-1})(s_{-1})\right] \\
					&= \sum_{s_1 \in A_1}\sum_{s_{-1} \in A_{-1}}\left[\sigma_1(s_1)f(\sigma_{-1})(s_{-1})\right] \\
					&= \sum_{s_1 \in A_1}\left[\sigma_1(s_1)\left(\sum_{s_{-1} \in A_{-1}}f(\sigma_{-1})(s_{-1})\right)\right] \\
					&= \left(\sum_{s_1 \in A_1}\sigma_1(s_1)\right)\left(\sum_{s_{-1} \in A_{-1}}f(\sigma_{-1})(s_{-1})\right) \\
					&= \sum_{s_{-1} \in A_{-1}}f(\sigma_{-1})(s_{-1}) \\
					&= \sum_{s_{-1} \in A_{-1}} \left[\sigma_2({(s_{-1})}_2)f(\sigma_{-1})(s_{-1})\right] \\
					&= \sum_{s_2 \in A_2}\sum_{s_{-1-2} \in A_{-1-2}} \left[\sigma_2(s_{2})f(\sigma_{-1-2})(s_{-1-2})\right] \\
					&= \left(\sum_{s_2 \in A_2} \sigma_2(s_{2})\right) \left(\sum_{s_{-1-2} \in A_{-1-2}} f(\sigma_{-1-2})(s_{-1-2})\right) \\
					&= \sum_{s_{-1-2} \in A_{-1-2}} f(\sigma_{-1-2})(s_{-1-2}) \\
					&\hspace{2mm} \vdots \\
					&= \sum_{s_{-1 - \ldots - (n-1)} \in A_{-1 - \ldots - (n-1)}} f(\sigma_{-1 - \ldots - (n-1)})(s_{-1 - \ldots - (n-1)}) \\
					&= \sum_{s_n \in A_n} \sigma_n(s_n) = 1.
				\end{align*}
				\item For each $\sigma, \sigma' \in \nabla(A)$ and $i \in N$, recalling that $\sigma_i(s_i) = 1 - \sum_{s'_i \in A_i-\{s_i\}}\sigma_i(s'_i)$ and $\sigma'_i(s_i) = 1 - \sum_{s'_i \in A_i-\{s_i\}}\sigma'_i(s'_i)$, we have:
				\begin{align*}
					f(\sigma) &= f(\sigma') \\
					\Rightarrow f(\sigma)(s) &= f(\sigma')(s) \text{ for all } s \in A \\
					\Rightarrow \prod_{i \in N}\sigma_i(s_i) &= \prod_{i \in N}\sigma'_i(s_i) \text{ for all } s \in A.
				\end{align*}
				If we fix $s_1 \in A_1$ then for all $s_{-1} \in A_{-1}$ we have:
				\[ \Rightarrow
				\begin{cases}
					\displaystyle\sigma_1(s'_1)f(\sigma_{-1})(s_{-1}) = \sigma'_1(s'_1)f(\sigma'_{-1})(s_{-1}) \text{ for all } s'_1 \in A_1-\{s_1\}& \\
					\displaystyle\left(1- \sum_{s'_1 \in A_1-\{s_1\}}\sigma_1(s'_1)\right)f(\sigma_{-1})(s_{-1}) = \left(1- \sum_{s'_1 \in A_1-\{s_1\}}\sigma'_1(s'_1)\right)f(\sigma'_{-1})(s_{-1})&
				\end{cases}
				\]
				\[ \Rightarrow
				\begin{cases} 
					\displaystyle\sum_{s'_1 \in A_1-\{s_1\}}\sigma_1(s'_1)f(\sigma_{-1})(s_{-1}) = \sum_{s'_1 \in A_1-\{s_1\}}\sigma'_1(s'_1)f(\sigma'_{-1})(s_{-1})& \\
					\displaystyle\sum_{s'_1 \in A_1-\{s_1\}}\sigma_1(s'_1)f(\sigma_{-1})(s_{-1}) = f(\sigma_{-1})(s_{-i}) - f(\sigma'_{-1})(s_{-1}) + \sum_{s'_1 \in A_1-\{s_1\}}\sigma'_1(s'_1)f(\sigma'_{-1})(s_{-1})&
				\end{cases}
				\]
				\begin{align*}
					\Rightarrow \sum_{s'_1 \in A_1-\{s_1\}}\sigma'_1(s'_1)f(\sigma'_{-1})(s_{-1}) &= f(\sigma_{-1})(s_{-1}) - f(\sigma'_{-1})(s_{-1}) + \sum_{s'_1 \in A_1-\{s_1\}}\sigma'_1(s'_1)f(\sigma'_{-1})(s_{-1}) \\
					\Rightarrow f(\sigma_{-1})(s_{-1}) &= f(\sigma'_{-1})(s_{-1}).
				\end{align*}
				Repeating the steps taken so far for all players excluding $i$ we get:
				\begin{align*}
					f(\sigma)(s) &= f(\sigma')(s) \text{ for all } s \in A\\
					\Rightarrow f(\sigma_{-1})(s_{-1}) &= f(\sigma'_{-1})(s_{-1}) \text{ for all } s_{-1} \in A_{-1} \\
					&\hspace{2mm} \vdots \\
					\Rightarrow f(\sigma_{-1-\ldots-(i-1)})(s_{-1-\ldots-(i-1)}) &= f(\sigma'_{-1-\ldots-(i-1)})(s_{-1-\ldots-(i-1)}) \\
					\text{ for all } s_{-1-\ldots-(i-1)} &\in A_{-1-\ldots-(i-1)} \\
					\Rightarrow f(\sigma_{-1-\ldots-(i-1)-(i+1)})(s_{-1-\ldots-(i-1)-(i+1)}) &= f(\sigma'_{-1-\ldots-(i-1)-(i+1)})(s_{-1-\ldots-(i-1)-(i+1)}) \\
					\text{ for all } s_{-1-\ldots-(i-1)-(i+1)} &\in A_{-1-\ldots-(i-1)-(i+1)} \\
					&\hspace{2mm} \vdots \\
					\hspace{-1cm}\Rightarrow f(\sigma_{-1-\ldots-(i-1)-(i+1)-\ldots-n})(s_{-1-\ldots-(i-1)-(i+1)-\ldots-n}) &= f(\sigma'_{-1-\ldots-(i-1)-(i+1)-\ldots-n})(s_{-1-\ldots-(i-1)-(i+1)-\ldots-n}) \\
					\text{ for all } s_{-1-\ldots-(i-1)-(i+1)-\ldots-n} &\in A_{-1-\ldots-(i-1)-(i+1)-\ldots-n} \\
					\Rightarrow \sigma_i(s_i) &= \sigma'_i(s_i) \text{ for all } s_i \in A_i \\
					\Rightarrow \sigma_i &= \sigma_i'.
				\end{align*}
				\item First let $\sigma \in \Delta(A)$. If we let $s \in A$ and specify values $\sigma(s') \geq 0$ for all $s' \in A-\{s\}$ such that $\sum_{s' \in A-\{s\}}\sigma(s') \leq 1$ then $\sigma$ is uniquely determined, as $\sigma(s) = 1 - \sum_{s' \in A-\{s\}}\sigma(s')$. Further, if we let $\sigma' \in \Delta(A)$ where there exists $s' \in A-\{s\}$ such that $\sigma'(s') \neq \sigma(s')$ then we trivially have that $\sigma' \neq \sigma$. Hence $\sigma$ is uniquely determined if and only if we have have values specified for all $s' \in A-\{s\}$, which is $-1 + |A| = -1 + \prod_{i \in N}|A_i|$ values. To complete the proof it suffices for us to show that fewer values will uniquely determine an arbitrary $\sigma \in \nabla(A)$.
				
				Let $\sigma = (\sigma_1, \ldots, \sigma_n) \in \nabla(A)$, $s \in A$ and $i \in N$. For each $j \in N$ pick $s'_j \in A_j-\{s_j\}$. Note that $\sigma_j(s'_j) = 1 - \sum_{s''_j \in A_j-\{s_j'\}}\sigma_j(s''_j)$. Now:
				\begin{enumerate}
					\item for each $s''_i \in A_i-\{s'_i\}$, suppose $f(\sigma)(s''_i, s_{-i})$ is specified; and 
					\item for each $j \in N_{-i}$ and $s''_j \in A_j$, suppose $f(\sigma)(s''_j, s_{-j})$ is specified.
				\end{enumerate}
				This gives us $\sum_{i \in N}\left(|A_i|-1\right) = -n + \sum_{i \in N}|A_i|$ values specified. Since $-n < -1$ for all $n \geq 2$ and $\sum_{i \in N}|A_i| \leq \prod_{i \in N}|A_i|$ for all $n$ and all $|A_1|, \ldots, |A_n| \geq 2$, it follows that $-n + \sum_{i \in N}|A_i| < -1 + \prod_{i \in N}|A_i|$ for all $n \geq 2$ and all $|A_1|, \ldots, |A_n| \geq 2$. 
				
				For each $j \in N_{-i}$ and $s''_j \in A_j$, $f(\sigma)(s''_j, s_{-j}) = \sigma_j(s''_j)\prod_{k \in N_{-j}}\sigma_k(s_k) = \sigma_i(s_i)\sigma_j(s''_j)\sum_{k \in N_{-i-j}}\sigma_k(s_k)$ with $f(\sigma)(s'_j, s_{-j}) = \sigma_i(s_i)\left(1 - \sum_{s''_j \in A_j-\{s_j'\}}\sigma_j(s''_j)\right)\sum_{k \in N_{-i-j}}\sigma_k(s_k)$. Rearranging we get:
				\begin{align}
					\sigma_i(s_i) &= \frac{f(\sigma)(s''_j, s_{-j})}{\sigma_j(s''_j)\prod_{k \in N_{-i-j}}\sigma_k(s_k)} \text{ for all } s''_j \in A_j; \text{ and} \label{eqn:zzz1} \\
					\sigma_i(s_i) &= \frac{f(\sigma)(s'_j, s_{-j})}{\left(1 - \sum_{s''_j \in A_j-\{s_j'\}}\sigma_j(s''_j)\right)\prod_{k \in N_{-i-j}}\sigma_k(s_k)}. \label{eqn:zzz2}
				\end{align}
				
				For each $s''_j \in A_j-\{s'_j\}$, if we set Equation \eqref{eqn:zzz1} equal to itself for the $s_j$ and $s_j''$ cases then rearrange we get:
				\begin{align}
					\frac{f(\sigma)(s_j, s_{-j})}{\sigma_j(s_j)\prod_{k \in N_{-i-j}}\sigma_k(s_k)} &= \frac{f(\sigma)(s''_j, s_{-j})}{\sigma_j(s''_j)\prod_{k \in N_{-i-j}}\sigma_k(s_k)} \nonumber\\
					\Rightarrow \sigma_j(s''_j) &= \frac{f(\sigma)(s''_j, s_{-j})}{f(\sigma)(s_j, s_{-j})}\sigma_j(s_j). \label{eqn:zzz3}
				\end{align}
				
				Subbing Equation \eqref{eqn:zzz3} in to Equation \eqref{eqn:zzz2} then setting equal to Equation \eqref{eqn:zzz1} for the $s_j$ case we get:
				\begin{align*}
					\frac{f(\sigma)(s'_j, s_{-j})}{\left(1 - \displaystyle\sum_{s''_j \in A_j-\{s'_j\}}\frac{f(\sigma)(s''_j, s_{-j})}{f(\sigma)(s_j, s_{-j})}\sigma_j(s_j)\right)\displaystyle\prod_{k \in N_{-i-j}}\sigma_k(s_k)} &= \frac{f(\sigma)(s_j, s_{-j})}{\sigma_j(s_j)\displaystyle\prod_{k \in N_{-i-j}}\sigma_k(s_k)} 
				\end{align*}
				\begin{align}
					\Rightarrow f(\sigma)(s'_j, s_{-j})\sigma_j(s_j) &= f(\sigma)(s_j, s_{-j}) - \sum_{s''_j \in A_j-\{s'_j\}}f(\sigma)(s''_j, s_{-j})\sigma_j(s_j) \nonumber \\
					\Rightarrow \sigma_j(s_j) &= \frac{f(\sigma)(s_j, s_{-j})}{\displaystyle\sum_{s''_j \in A_j}f(\sigma)(s''_j, s_{-j})} = \frac{f(\sigma)(s_j, s_{-j})}{\displaystyle\sum_{s'''_j \in A_j}f(\sigma)(s'''_j, s_{-j})}. \label{eqn:zzz4}
				\end{align}
				Subbing Equation \eqref{eqn:zzz4} in to Equation \eqref{eqn:zzz3}, for each $s''_j \in A_j-\{s'_j\}$ we get:
				\[
					\sigma_j(s''_j) = \frac{f(\sigma)(s''_j, s_{-j})}{f(\sigma)(s_j, s_{-j})}\sigma_j(s_j) = \frac{f(\sigma)(s''_j, s_{-j})}{f(\sigma)(s_j, s_{-j})}\frac{f(\sigma)(s_j, s_{-j})}{\displaystyle\sum_{s'''_j \in A_j}f(\sigma)(s'''_j, s_{-j})} = \frac{f(\sigma)(s''_j, s_{-j})}{\displaystyle\sum_{s'''_j \in A_j}f(\sigma)(s'''_j, s_{-j})}.
				\]
				Which gives us:
				\begin{align*}
					\sigma_j(s'_j) &= 1 - \sum_{s''_j \in A_j-\{s'_j\}}\sigma_j(s''_j) \\
					&= 1 - \sum_{s''_j \in A_j-\{s'_j\}}\frac{f(\sigma)(s''_j, s_{-j})}{\displaystyle\sum_{s'''_j \in A_j}f(\sigma)(s'''_j, s_{-j})}\\
					&= \frac{\left(\displaystyle\sum_{s'''_j \in A_j}f(\sigma)(s'''_j, s_{-j})\right) - \left(\displaystyle\sum_{s''_j \in A_j-\{s'_j\}}f(\sigma)(s''_j, s_{-j})\right)}{\displaystyle\sum_{s'''_j \in A_j}f(\sigma)(s'''_j, s_{-j})} \\
					&= \frac{f(\sigma)(s'_j, s_{-j})}{\displaystyle\sum_{s'''_j \in A_j}f(\sigma)(s'''_j, s_{-j})}. 
				\end{align*}
				
				Finally, for each $s''_i \in A_i-\{s'_i\}$ we have:
				\begin{align*}
					\sigma_i(s''_i) &= \frac{f(\sigma)(s''_i, s_{-i})}{\displaystyle\prod_{j \in N_{-i}}\sigma_j(s_j)}; \text{ and} \\
					\sigma_i(s'_i) &= 1 - \sum_{s''_i \in A_i-\{s'_i\}}\sigma_i(s''_i) = 1 - \frac{\displaystyle\sum_{s''_i \in A_i-\{s'_i\}}f(\sigma)(s''_i, s_{-i})}{\displaystyle\prod_{j \in N_{-i}}\sigma_j(s_j)}.
				\end{align*}
			\end{enumerate}
			
			Hence with $-n + \sum_{i \in N}|A_i|$ values already specified the rest of $\sigma$ is uniquely determined. Consequently $\Delta(A) - f(\nabla(A))$ is non-empty. 
		\end{proof}
	\end{proposition}
	
	To simplify notation for the remainder of the paper, for each $\sigma \in \nabla(A)$ and $s \in A$ we shall denote $f(\sigma)(s)$ as $\sigma(s)$. 
	
	For each player $i \in N$, the domain for their utility function can be extended linearly from $A$ to $\nabla(A)$ with $\textswab{u}_i(\sigma) = \Sigma_{s \in A}\sigma(s)u_i(s)$ for all $\sigma \in \nabla(A)$.
	
	A \textit{pure strategy Nash equilibrium} is a strategy profile $s \in A$ where for each $i \in N$, $u_i(s_i, s_{-i}) \geq u_i(s_i', s_{-i})$ for all $s_i' \in A_i$. For example, in Example \ref{fullsymeg} the profile $(b,b,b)$ is a pure strategy Nash equilibrium.
	
	For each player $i \in N$, the \textit{maximin value} for player $i$ is given by:
	\[\underline{u}_i = \max_{\sigma_i \in \Delta(A_i)}\min_{\sigma_{-i} \in {\nabla(A)}_{-i}} u_i(\sigma_i, \sigma_{-i}),\]
	and the \textit{minimax value} for player $i$ is given by:
	\[\overline{u}_i = \min_{\sigma_{-i} \in {\nabla(A)}_{-i}}\max_{\sigma_i \in \Delta(A_i)} u_i(\sigma_i, \sigma_{-i}).\]
	
	The maximin value for player $i$ from an intuitive point of view is:
	\begin{enumerate}
		\item the highest expected payoff player $i$ can be sure to obtain when they do not know what strategies their opponents will play; and
		\item the lowest expected payoff player $i$'s opponents can force player $i$ to obtain when they know what strategy player $i$ will play.
	\end{enumerate}
	
	While the minimax value for player $i$ from an intuitive point of view is:
	\begin{enumerate}
		\item the lowest expected payoff player $i$'s opponents can force player $i$ to obtain when they do not know what strategy player $i$ will play; and
		\item the highest expected payoff player $i$ can be sure to obtain when they know what strategies their opponents will play.
	\end{enumerate} 
	
	It is obvious from both intuitive points of view for the maximin and minimax values that the maximin values are less than or equal to the minimax values, ie. for each player $i \in N$:
	\[\underline{u}_i = \max_{\sigma_i \in \Delta(A_i)}\min_{\sigma_{-i} \in {\nabla(A)}_{-i}} u_i(\sigma_i, \sigma_{-i}) \leq \min_{\sigma_{-i} \in {\nabla(A)}_{-i}}\max_{\sigma_i \in \Delta(A_i)} u_i(\sigma_i, \sigma_{-i}) = \overline{u}_i.\]
	
	A fairly standard proof of the above is as follows.
	
	\begin{proposition}
		Given a game $\Gamma = (N, A, u)$, for each player $i \in N$:
		\[\underline{u}_i = \max_{\sigma_i \in \Delta(A_i)}\min_{\sigma_{-i} \in {\nabla(A)}_{-i}} u_i(\sigma_i, \sigma_{-i}) \leq \min_{\sigma_{-i} \in {\nabla(A)}_{-i}}\max_{\sigma_i \in \Delta(A_i)} u_i(\sigma_i, \sigma_{-i}) = \overline{u}_i.\]
		
		\begin{proof}
			Let $\sigma^*_i = \displaystyle\argmax_{\sigma_i \in \Delta(A_i)}\min_{\sigma_{-i} \in {\nabla(A)}_{-i}} u_i(\sigma_i, \sigma_{-i})$ and $\sigma^*_{-i} = \displaystyle\argmin_{\sigma_{-i} \in {\nabla(A)}_{-i}}\max_{\sigma_i \in \Delta(A_i)} u_i(\sigma_i, \sigma_{-i})$.
			
			\begin{enumerate}
				\item $\underline{u}_i = \displaystyle\max_{\sigma_i \in \Delta(A_i)}\min_{\sigma_{-i} \in {\nabla(A)}_{-i}} u_i(\sigma_i, \sigma_{-i}) \leq u_i(\sigma^*_i, \sigma_{-i})$ for all $\sigma_{-i} \in {\nabla(A)}_{-i}$; and
				\item $u_i(\sigma_i, \sigma^*_{-i}) \leq \displaystyle\min_{\sigma_{-i} \in {\nabla(A)}_{-i}}\max_{\sigma_i \in \Delta(A_i)} u_i(\sigma_i, \sigma_{-i}) = \overline{u}_i$ for all $\sigma_i \in \Delta(A_i)$.
			\end{enumerate}
			
			From which it follows that $\underline{u}_i \leq u_i(\sigma^*_i, \sigma^*_{-i}) \leq \overline{u}_i$.
		\end{proof}
	\end{proposition}
	
	A game $\Gamma = (N, A, u)$ is referred to as \textit{zero-sum} if $\sum_{i \in N} u_i(s) = 0$ for all $s \in A$. In a zero-sum game, each player's gain or loss is exactly matched by the combined losses or gains of their opponents. 
	
	\begin{proposition}
		Given a $2$-player zero-sum game $(\{1, 2\}, A_1 \times A_2, (u_1, u_2))$, the maximin and minimax values for each player are equal to minus the minimax and maximin values for the other player.
		
		\begin{proof}
			The following is rephrased from various different sources. This may be seen by noting that $u_1(s) = -u_2(s)$, and $u_2(s) = -u_1(s)$, for all $s \in A$. Hence:
			\begin{align*}
				\text{(i) } \underline{u}_1 = \max_{\sigma_1 \in \Delta(A_1)}\min_{\sigma_2 \in \Delta(A_2)} u_1(\sigma_1, \sigma_2) &= \max_{\sigma_1 \in \Delta(A_1)}\min_{\sigma_2 \in \Delta(A_2)} -u_2(\sigma_1, \sigma_2) \\
				&= -\min_{\sigma_1 \in \Delta(A_1)}\max_{\sigma_2 \in \Delta(A_2)} u_2(\sigma_1, \sigma_2) = -\overline{u}_2; \text{ and} \\
				\text{(ii) } \overline{u}_1 = \min_{\sigma_1 \in \Delta(A_1)}\max_{\sigma_2 \in \Delta(A_2)} u_1(\sigma_1, \sigma_2) &= \min_{\sigma_1 \in \Delta(A_1)}\max_{\sigma_2 \in \Delta(A_2)} -u_2(\sigma_1, \sigma_2) \\
				&= -\max_{\sigma_1 \in \Delta(A_1)}\min_{\sigma_2 \in \Delta(A_2)} u_2(\sigma_1, \sigma_2) = -\underline{u}_2. 
			\end{align*}
		\end{proof}
	\end{proposition}
	
	\begin{proposition} \label{prop:2pzsminmax=maxmin}
		\cite{VNM} Given a $2$-player zero-sum game $\Gamma = (\{1, 2\}, A_1\times A_2, (u_1, u_2))$, the maximin and minimax values for each player are equal. Ie. $\underline{u}_i = \overline{u}_i$ for all $i \in \{1, 2\}$.
	\end{proposition}
	
	We denote the subgroup relation as $\leq$, the group generated by a subset $H$ of a group $G$ as $\langle{H}\rangle$, the group of permutations on a non-empty set $X$ as $S_X$, and the subset of transpositions on $X$ as $T_X$. The reader is reminded that the permutations on $X$ are equivalent to the bijections from $X$ to itself, henceforth we will refer to them interchangeably.
	
	An \textit{action} of a group $G$ on a set $N$ is a homomorphism $\alpha$ from $G$ into the bijections from $N$ to itself. For each $g \in G$ and $i \in N$ we denote $\bigl(\alpha(g)\bigr)(i)$ as $g(i)$. When $G$ acts on the left or right of $N$ the action is called a left or right action respectively. We note that left actions can be defined equivalently as antihomomorphisms that act on the right, and dually for right actions.
	
	An action is \textit{transitive} if for each $i, j \in N$ there exists $g \in G$ such that $g(i) = j$, \textit{regular} if for each $i, j \in N$ there exists precisely one $g \in G$ such that $g(i) = j$, and \textit{$n$-transitive} if for each $\pi \in S_N$ there exists $g \in G$ such that $g(i) = \pi(i)$ for all $i \in N$. When an action of $G$ can be inferred we simply refer to $G$ as being transitive, regular or $n$-transitive respectively.
	
	The \textit{stabiliser of $i \in N$}, which we denote as $G_i$, is the subgroup $\{g \in G: g(i) = i\}$ of elements in $G$ that fix $i$. Similarly the \textit{stabiliser of $N$}, which we denote as $G_N$, is the normal subgroup $\{g \in G: g(i) = i \text{ for all } i \in N\} = \cap_{i \in N}G_i$ of elements in $G$ that fix each $i \in N$.
	
	The \textit{orbit of $i \in N$} is $G(i) = \{g(i): g \in G\}$. The \textit{orbits of $N$}, denoted as $N/G$, is the set $\{G(i): i \in N\}$ which forms a partition of $N$.
	
	By a \textit{groupoid} we mean a category in which every morphism is invertible. For the sake of brevity, when the objects of a groupoid can be inferred we refer to the morphisms as a groupoid. 
    \section{Label-Dependent Notions of Symmetry} \label{sec:labdepnotions}
There are various ways to define a notion of symmetry, not all of which are distinct. In each case we need all players to have the same number of strategies, consequently all games are implicitly $m$-strategy games. It is often assumed when defining symmetric games that all players have the same strategy labels and any notion of symmetry will treat the same labels as equivalent. We shall refer to these as \textit{label-dependent} notions. 

\subsection{Permutations Acting On Strategy Profiles}
There is some confusion over how to correctly define symmetric games, see \cite[Definition 7]{DMaskin}, in order to provide clarity we need to review two ways that player permutations may act on strategy profiles. 

Given a player permutation $\pi \in S_N$ and strategy profile $s \in A$, two possible action choices are $(s_i)_{i \in N} \mapsto (s_{\pi(i)})_{i \in N}$ and $(s_i)_{i \in N} \mapsto (s_{\pi^{-1}(i)})_{i \in N}$. We denote $(s_{\pi^{-1}(i)})_{i \in N}$ as $\pi(s)$, for example given $(s_1, \ldots, s_n) \in A$, $\pi(s_1, \ldots, s_n) = (s_{\pi^{-1}(1)}, \ldots, s_{\pi^{-1}(n)})$.

The author notes that our somewhat unintuitive notation has been chosen so that it matches with composition and inversion in an ideal manner. That is so for each $s \in A$, $(\tau \circ \pi)(s) = \tau\bigl(\pi(s)\bigr)$ and $(\tau \circ \pi)^{-1} = \pi^{-1} \circ \tau^{-1}$. 

\begin{lemma} \label{simpleactionprop}
	For each $s \in A$ and $\pi \in S_N$, $s \mapsto \pi(s)$ is a left action of $S_N$ on $A$.
	\begin{proof}
		The identity permutation trivially acts as an identity so we need only establish associativity. For each $\pi, \tau \in S_N$, $s \in A$ and $i \in N$, $\bigl((\tau \circ \pi)(s)\bigr)_i = s_{(\tau \circ \pi)^{-1}(i)} = s_{\pi^{-1}(\tau^{-1}(i))} = \bigl(\pi(s)\bigr)_{\tau^{-1}(i)} = \Bigl(\tau\bigl(\pi(s)\bigr)\Bigr)_i$.
		
		For each $s, s' \in A$. We have $\pi(s) = \pi(s')$ if and only if $s_{\pi^{-1}(i)} = s_{\pi^{-1}(i)}'$ for all $i \in N$ and for each $s \in A$, $\pi^{-1}(s) \in A$ and $\pi(\pi^{-1}(s)) = (\pi \circ \pi^{-1})(s) = s$. Hence for each $s \in A$, $s \mapsto \pi(s)$ is both injective and surjective, ie. $s \mapsto \pi(s) \in \bij(A, A)$.
	\end{proof}
\end{lemma}

It might be worth explaining which bijections of $A$ we have using $\{s\mapsto\pi(s): \pi \in S_N\}$? Ie. $\bij(A)-\{s\mapsto\pi(s): \pi \in S_N\}$, though is trivially just the bijections where there's consistency with mapping strategies from one player to the same possibly other player. Same for $\{s\mapsto g(s): g \in \bij(\Gamma)\}$ later. 

Since $\pi^{-1}(s) = (s_{\pi(i)})_{i \in N}$ for all $s \in A$, $s \mapsto \pi(s)$ and $s \mapsto \pi^{-1}(s)$ are dual to each other. Hence the dual results hold for $\pi^{-1}$.

\begin{lemma} 
	For each $s \in A$ and $\pi \in S_N$, $s \mapsto \pi^{-1}(s)$ is a right action of $S_N$ on $A$.
\end{lemma}

Given $\pi \in S_N$ we denote the map $s \mapsto u_{\pi(i)}\bigl(\pi(s)\bigr)$ as $u_{\pi(i)} \circ \pi$. Note that $u_{\pi(i)} \circ \pi$ is the utility function of player $\pi(i)$ when the strategy profiles are acted upon by the player permutation $\pi$.

\begin{corollary} \label{utilityactionprop}
	For each $\pi, \tau \in S_N$, $u_{(\tau \circ \pi)(i)} \circ (\tau \circ \pi) = (u_{\tau(\pi(i))} \circ \tau) \circ \pi$.
	\begin{proof}
		For each $i \in N$, $s \in A$, 
		\begin{align*}
			\left(u_{(\tau \circ \pi)(i)} \circ (\tau \circ \pi)\right)(s) &= u_{(\tau \circ \pi)(i)}\Bigl((\tau \circ \pi)(s)\Bigr) \\
			&= u_{\tau(\pi(i))}\Bigl(\tau\bigl(\pi(s)\bigr)\Bigr) \\
			&= \left((u_{\tau(\pi(i))} \circ \tau) \circ \pi\right)(s).
		\end{align*}
	\end{proof}
\end{corollary}

The above may all be done exactly the same but with mixed strategy profiles, however the bulk majority of cases where we use permutations acting on strategy profiles it will be for pure strategy profiles. Oppositely, the following may all be done exactly the same but with pure strategy profiles, however the purposes that led to bothering with such notation in the first place involves mixed not pure strategy profiles. 

Given $\pi \in S_N$ and $\sigma \in \nabla(A)$, we denote:
\begin{enumerate}
	\item $\sigma_i$ as $\pi(\sigma_i)$; and
	\item $(\sigma_{\pi^{-1}(1)}, \ldots, \sigma_{\pi^{-1}\left(\pi(i)-1\right)}, \sigma_{\pi^{-1}\left(\pi(i)+1\right)}, \ldots, \sigma_{\pi^{-1}(n)})$ as $\pi(\sigma_{-i})$.
\end{enumerate}
This gives us $\pi(\sigma_i) = \pi(\sigma)_{\pi(i)} \in \Delta(A_{\pi(i)})$ and $\pi(\sigma_{-i}) = \pi(\sigma)_{-\pi(i)} \in {\nabla(A)}_{-\pi(i)}$, however $\pi(\sigma_i)$ does not need a $\sigma_{-i} \in {\nabla(A)}_{-i}$ unlike $\pi(\sigma)_{\pi(i)}$, similarly $\pi(\sigma_{-i})$ does not need a $\sigma_i \in \Delta(A_i)$ unlike $\pi(\sigma)_{-\pi(i)}$, which will make the proof of Proposition \ref{prop:standsymgamesaremaxminandminmaxfair} much easier to understand.

\begin{proposition}
	For each $i \in N$, $\pi \in S_N$ and $\sigma \in \nabla(A)$, $\pi(\sigma_i, \sigma_{-i}) = \Bigl(\pi(\sigma_i), \pi(\sigma_{-i})\Bigr)$.
	
	\begin{proof}
		Note that for each $(\sigma_i, \sigma_{-i}) \in \Delta(A_i)\times{\nabla(A)}_{-i}$, $(\sigma_i, \sigma_{-i}) = \sigma = (\sigma_j, \sigma_{-j}) \in \Delta(A_j)\times{\nabla(A)}_{-j}$.
		
		Now, $\pi(\sigma_i, \sigma_{-i}) = \pi(\sigma) = \Bigl(\pi(\sigma)_{\pi(i)}, \pi(\sigma)_{-\pi(i)}\Bigr) = \Bigl(\pi(\sigma_i), \pi(\sigma_{-i})\Bigr)$.
	\end{proof}
\end{proposition}

\begin{proposition} \label{prop:firstgroupoidprop}
	$\set{\sigma_i \mapsto \pi(\sigma_i): i \in N, \pi \in S_N}$ is a subgroupoid of $\cup_{i, j \in N}\bij(\Delta(A_i), \Delta(A_j))$.
	
	\begin{proof}		
		Let $Y = \set{\sigma_i \mapsto \pi(\sigma_i): i \in N, \pi \in S_N}$.
		
		\begin{enumerate}
			\item For each $i \in N$, $\id_N(\sigma_i) = \sigma_i$ for all $\sigma_i \in \Delta(A_i)$. Hence $\id_{\Delta(A_i)} = \sigma_i \mapsto \id_N(\sigma_i) \in Y$;
			\item For each $i \in N$ and $\pi \in S_N$, we trivially have $\sigma_i \mapsto \pi(\sigma_i) \in \bij(\Delta(A_i), \Delta(A_{\pi(i)}))$;
			\item For each $i \in N$ and $\pi, \tau \in S_N$, $\tau(\pi(\sigma_i)) = \tau(\sigma_i) = \sigma_i = (\tau \circ \pi)(\sigma_i)$. Hence $(\sigma_{\pi(i)} \mapsto \tau(\sigma_{\pi(i)})) \circ (\sigma_i \mapsto \pi(\sigma_i)) = \sigma_i \mapsto (\tau \circ \pi)(\sigma_i) \in Y$; and
			\item Finally, for each $i \in N$ and $\pi \in S_N$, since $\pi \circ \pi^{-1} = \id_N = \pi^{-1} \circ \pi$, we have $(\sigma_i \mapsto \pi(\sigma_i))^{-1} = \sigma_{\pi(i)} \mapsto \pi^{-1}(\sigma_{\pi(i)}) \in Y$.
		\end{enumerate} \vspace{-0.8cm}
	\end{proof}
\end{proposition}

\begin{proposition} \label{prop:secondgroupoidprop}
	$\set{\sigma_{-i} \mapsto \pi(\sigma_{-i}): i \in N, \pi \in S_N}$ is a subgroupoid of $\cup_{i, j \in N}\bij({\nabla(A)}_{-i}, {\nabla(A)}_{-j})$.
	
	\begin{proof}
		Let $Y = \set{\sigma_{-i} \mapsto \pi(\sigma_{-i}): i \in N, \pi \in S_N}$.

		\begin{enumerate}
			\item For each $i \in N$, it follows trivially from $\id_N^{-1} = \id_N$ that for each $\sigma_{-i} \in {\nabla(A)}_{-i}$:
			\begin{align*}
				\textstyle\id_N(\sigma_{-i}) &= (\sigma_{\id_N^{-1}(1)}, \ldots, \sigma_{\id_N^{-1}\left(\id_N(i)-1\right)}, \sigma_{\id_N^{-1}\left(\id_N(i)+1\right)}, \ldots, \sigma_{\id_N^{-1}(n)})\\
				&= (\sigma_1, \ldots, \sigma_{\id_N(i)-1}, \sigma_{\id_N(i)+1}, \ldots, \sigma_n)\\
				&= (\sigma_1, \ldots, \sigma_{i-1}, \sigma_{i+1}, \ldots, \sigma_n) \\
				&= \sigma_{-i}.
			\end{align*}
			Hence $\id_{{\nabla(A)}_{-i}} = \sigma_{-i} \mapsto \id_N(\sigma_{-i}) \in Y$;
			\item Let $i \in N$, $\pi \in S_N$ and $\sigma_{-i}, \sigma'_{-i} \in {\nabla(A)}_{-i}$, then:
			\begin{align*}
				\pi(\sigma_{-i}) &= \pi(\sigma'_{-i}) \\
				\Rightarrow \sigma_{\pi^{-1}(j)} &= \sigma'_{\pi^{-1}(j)} \text{ for all } j \in \set{1, \ldots, \pi(i)-1, \pi(i)+1, \ldots, n} \\
				\Rightarrow \sigma_j &= \sigma'_j \text{ for all } j \in \set{1, \ldots, i-1, i+1, \ldots, n}.
			\end{align*}			
			Therefore $\sigma_{-i} \mapsto \pi(\sigma_{-i})$ is injective. Now let $i \in N$, $\sigma_{-\pi(i)} = (\sigma_1, \ldots, \sigma_{\pi(i)-1}, \sigma_{\pi(i)+1}, \ldots, \sigma_n) \in {\nabla(A)}_{-\pi(i)}$. Note $(\sigma_{\pi(1)}, \ldots, \sigma_{\pi(i)-1}, \sigma_{\pi(i)+1}, \ldots, \sigma_{\pi(n)}) \in {\nabla(A)}_{-i}$ and:
			\begin{align*}
				\pi(\sigma_{\pi(1)}, \ldots, \sigma_{\pi(i)-1}, \sigma_{\pi(i)+1}, \ldots, \sigma_{\pi(n)}) &= (\sigma_{\pi^{-1}(\pi(1))}, \ldots, \sigma_{\pi^{-1}(\pi(\pi(i)-1))}, \sigma_{\pi^{-1}(\pi(\pi(i)+1))}, \ldots, \sigma_{\pi^{-1}(\pi(n))})\\
				&= (\sigma_1, \ldots, \sigma_{\pi(i)-1}, \sigma_{\pi(i)+1}, \ldots, \sigma_n)\\
				&= \sigma_{-\pi(i)}.
			\end{align*}
			It follows that $\sigma_{-i} \mapsto \pi(\sigma_{-i})$ is surjective, hence $\sigma_{-i} \mapsto \pi(\sigma_{-i}) \in \bij({\nabla(A)}_{-i}, {\nabla(A)}_{-\pi(i)})$;
			\item For each $i \in N$ and $\pi, \tau \in S_N$:
			\begin{align*}
				\tau(\pi(\sigma_{-i})) &= \tau(\sigma_{\pi^{-1}(1)}, \ldots, \sigma_{\pi^{-1}\left(\pi(i)-1\right)}, \sigma_{\pi^{-1}\left(\pi(i)+1\right)}, \ldots, \sigma_{\pi^{-1}(n)}) \\
				&= (\sigma_{\pi^{-1}(\tau^{-1}(1))}, \ldots, \sigma_{\pi^{-1}\left(\tau^{-1}(\tau(\pi(i))-1)\right)}, \sigma_{\pi^{-1}\left(\tau^{-1}(\tau(\pi(i))+1)\right)}, \ldots, \sigma_{\pi^{-1}(\tau^{-1}(n))}) \\
				&= (\sigma_{(\tau \circ \pi)^{-1}(1)}, \ldots, \sigma_{(\tau \circ \pi)^{-1}((\tau \circ \pi)(i))-1)}, \sigma_{(\tau \circ \pi)^{-1}((\tau \circ \pi)(i))+1)}, \ldots, \sigma_{(\tau \circ \pi)^{-1}(n))}) \\
				&= (\tau \circ \pi)(\sigma_{-i}).
			\end{align*}
			Hence $(\sigma_{-\pi(i)} \mapsto \tau(\sigma_{-\pi(i)})) \circ (\sigma_{-i} \mapsto \pi(\sigma_{-i})) = \sigma_{-i} \mapsto (\tau \circ \pi)(\sigma_{-i}) \in Y$; and
			\item Finally, for each $i \in N$ and $\pi \in S_N$, since $\pi \circ \pi^{-1} = \id_N = \pi^{-1} \circ \pi$, we have $(\sigma_{-i} \mapsto \pi(\sigma_{-i}))^{-1} = \sigma_{-\pi(i)} \mapsto \pi^{-1}(\sigma_{-\pi(i)}) \in Y$.
		\end{enumerate} \vspace{-0.8cm}
	\end{proof}
\end{proposition}

\subsection{Game Invariants}
Game invariants give us a notion of players being indifferent between the current positions and an alternative arrangement of positions.

\begin{definition}
	$\pi \in S_N$ is an \textit{invariant} of $\Gamma$ if for each $i \in N$, $u_i = u_{\pi(i)} \circ \pi$.
\end{definition}

\begin{lemma}
	The invariants of a game form a group.
	\begin{proof}
		Since the identity permutation $e \in S_N$ acts as an identity on $A$ it follows that $u_i = u_i \circ e$ for all $i \in N$, hence $e$ is an invariant. Suppose $\pi \in S_N$ is an invariant of $\Gamma$, and hence that for each $i \in N$, $u_{\pi^{-1}(i)} = u_i \circ \pi$. Then for each $i \in N$, $u_i = (u_i \circ \pi) \circ \pi^{-1} = u_{\pi^{-1}(i)} \circ \pi^{-1}$. Finally suppose $\pi, \tau \in S_N$ are invariants of $\Gamma$. Then for each $i \in N$, $u_i = u_{\pi(i)} \circ \pi = (u_{\tau(\pi(i))} \circ \tau) \circ \pi = u_{(\tau \circ \pi)(i)} \circ (\tau \circ \pi)$.
	\end{proof}
\end{lemma}

\subsection{Notions of Anonymity}
Before surveying label-dependent notions of fairness we review various notions of \textit{anonymity}, which have previously been examined by Brandt et al. \cite{brandt2009symmetries}. 

Central to anonymity is the notion that players do not distinguish between their opponents, by which we mean each player merely cares about the strategies being played by their opponents and is indifferent between who is playing them.

\begin{definition}
	$\Gamma$ is \textit{weakly anonymous} if for each $i \in N$, $\pi \in S_{N-\{i\}}$, $u_i = u_i \circ \pi$.
\end{definition}

\begin{example} \label{weaklyanoneg}
	Weakly Anonymous 3-player game.
	\begin{center}
	\begin{game}{2}{2}[$(a,,)$]
     \> $a$      \> $b$ \\
		$a$   \> $0,1,2$  \> $4,6,7$ \\
		$b$   \> $4,5,8$  \> $9,12,14$
	\end{game}
	\hspace*{10mm} 
	\begin{game}{2}{2}[$(b,,)$]
     \> $a$     \> $b$ \\
		$a$   \> $3,6,8$ \> $10,11,14$ \\
		$b$   \> $10,12,13$ \> $15,16,17$
	\end{game}
	\end{center}
	The reader may like to verify that $u_i = u_i \circ (jk)$ for all distinct $i, j, k \in N$. For example, $u_1(a,b,a) = 4 = u_1\bigl((23)(a,b,a)\bigr) = u_1(a,a,b)$. Since $S_{N-\{i\}} = \{e, (jk)\}$ for all $i \in N$, $\Gamma$ is weakly anonymous.
	
	When we say that players do not distinguish between their opponents, we mean for example that when playing $a$, player $1$ is indifferent between the strategy profiles $(a,a,b)$ and $(a,b,a)$.
\end{example}
	
Weak anonymity may be strengthened by requiring the players care merely about the strategies being played and be indifferent between who is playing each strategy, or equivalently, by requiring each player have the same payoff for each orbit in $A/S_N$.

\begin{definition}
	$\Gamma$ is \textit{anonymous} if for each $i \in N$, $\pi \in S_N$, $u_i = u_i \circ \pi$. 
\end{definition}

\begin{example} \label{anoneg}
		Anonymous 3-player game.
		\begin{center}
		\begin{game}{2}{2}[$(a,,)$]
			      \> $a$      \> $b$ \\
			$a$   \> $0,1,2$  \> $3,4,5$ \\
			$b$   \> $3,4,5$  \> $6,7,8$
		\end{game}
		\hspace*{10mm} 
		\begin{game}{2}{2}[$(b,,)$]
			      \> $a$     \> $b$ \\
			$a$   \> $3,4,5$ \> $6,7,8$ \\
			$b$   \> $6,7,8$ \> $9,10,11$
		\end{game}
		\end{center}
		The reader may like to verify the orbits of $A$ are given by $A/S_N = \bigl\{\{(a,a,a)\}, \newline\{(a,a,b), (a,b,a), (b,a,a)\}, \{(a,b,b), (b,a,b), (b,b,a)\}, \{(b,b,b)\}\bigr\}$ and that each player has the same payoff for each orbit in $A/S_N$. 
		
		For example, let $\pi = (123)$, then we have $\pi(s_1, s_2, s_3) = (s_{\pi^{-1}(1)}, s_{\pi^{-1}(2)}, s_{\pi^{-1}(3)}) = (s_3, s_1, s_2)$ giving us $\pi(a,a,b) = (b,a,a)$.
	\end{example}
	
Anonymity may be strengthened also by requiring all players have the same payoff for each orbit in $A/S_N$.
	
\begin{definition}
	$\Gamma$ is \textit{fully anonymous} if for each $i, j \in N$, $\pi \in S_N$, $u_i = u_j \circ \pi$. 
\end{definition}

\begin{example}
		Fully anonymous 3-player game.
		\begin{center}
		\begin{game}{2}{2}[$(a,,)$]
			      \> $a$      \> $b$ \\
			$a$   \> $1,1,1$  \> $2,2,2$ \\
			$b$   \> $2,2,2$  \> $3,3,3$
		\end{game}
		\hspace*{10mm} 
		\begin{game}{2}{2}[$(b,,)$]
			      \> $a$     \> $b$ \\
			$a$   \> $2,2,2$ \> $3,3,3$ \\
			$b$   \> $3,3,3$ \> $4,4,4$
		\end{game}
		\end{center}
		The orbits of $A$ for the above game are the same as in Example \ref{anoneg}, however now all players have the same payoff for each orbit.
	\end{example}
	
In a fully anonymous game each player is indifferent between which position they play. Hence fully anonymous games are one class of games that fall under fairness. 

Note that the published version of \cite{brandt2009symmetries} refers to weakly anonymous, anonymous and fully anonymous games as \textit{weakly symmetric}, \textit{weakly anonymous} and \textit{strongly anonymous games} respectively. The reason for this is the author finds using the symmetric terminology in the context of anonymity rather confusing when it is already the convention to use the term symmetric for notions of symmetry/fairness. Further, since the anonymity notion that \cite{brandt2009symmetries} refer to as weakly symmetric is not the weakest of the three anonymity notions, the author has instead chosen to refer to them as simply anonymous, and so what \cite{brandt2009symmetries} refer to as weakly symmetric the author refers to as weakly anonymous. The author uses fully anonymous instead of strongly anonymous to be consistent with the terminology used for notions of symmetry/fairness.

\subsection{Notions of Symmetry} \label{subsec:labeldepnotionsofsymmetry}
Our broad requirements for fairness that players be indifferent between which position they play may be made more precise by requiring the invariants of a game be a transitive subgroup of $S_N$.

\begin{definition}
	$\Gamma$ is \textit{standard symmetric} \cite{NoahXE} if there exists a transitive subgroup $H$ of the player permutations such that for each $i \in N$ and $\pi \in H$, $u_i = u_{\pi(i)} \circ \pi$. 
\end{definition}

	In a standard symmetric game, while being indifferent between which position they play, each player may care about the arrangement of their opponents, or alternatively may distinguish between their opponents.

	\begin{example} \label{stdsymeg} Standard symmetric 3-player game.
		\begin{center}
  		\begin{game}{2}{2}[$(a,,)$]
				\>  $a$      \>  $b$      \\
			$a$	\>  $1,1,1$  \>  $3,7,4$  \\
			$b$	\>  $7,4,3$  \>  $6,5,8$  
		\end{game}
		\hspace*{5mm}
		\begin{game}{2}{2}[$(b,,)$]
				\>  $a$      \>  $b$      \\
			$a$	\>  $4,3,7$  \>  $8,6,5$  \\
   			$b$	\>  $5,8,6$  \>  $2,2,2$  
		\end{game}
		\end{center}
		
		The reader may like to verify that $\Gamma$ is invariant under $(123)$ and not invariant under $(12)$. Since $\langle (123)\rangle = \{e, (123), (132)\}$ is a transitive subgroup of $S_3$, $\Gamma$ is standard symmetric. Furthermore since $(12)$ is not an invariant the players are not indifferent between all possible position arrangements.
		
		A useful analogy for considering the fairness of $\Gamma$ is a game with three players sitting in a circle such that each player is indifferent between circular rotations of positions, and not indifferent to their opponents swapping positions. A similar notion of fairness/symmetry is often used by the author when coding map generators for artificial intelligence programming contests where users write bots to play games against each other. The maps are two dimensional grids with the edges wrapped, ie. on the surface of a torus, and constructed in such a way that everyone is indifferent between some reorderings of the players.
	\end{example}
	
	We obtain a stronger level of fairness by requiring the players be indifferent between all possible position rearrangements, that is by requiring all player permutations be invariants.

\begin{definition} \label{fullsymdef}
	$\Gamma$ is \textit{fully symmetric} if it is invariant under $S_N$.
\end{definition}

The reader may like to verify that Example \ref{fullsymeg} is invariant under the permutations $(12)$ and $(123)$. For example, let $\pi = (123)$, then $\pi(s_1, s_2, s_3) = (s_3, s_1, s_2)$ giving us $u_1(b, a, a) = u_2(a, b, a) = u_3(a, a, b) = 3$. Since invariants are closed under composition and $\langle (12), (123)\rangle = S_3$, Example \ref{fullsymeg} is fully symmetric. 

Next we establish that Definition \ref{fullsymdef} can be characterised by various conditions. 

\begin{theorem} \label{basicsymequivthm}
	The following conditions are equivalent:
	\begin{enumerate}
		\item $\Gamma$ is fully symmetric;
		\item $\Gamma$ is standard symmetric and weakly anonymous;
		\item For each $i \in N$ and $\pi \in S_N$, $u_{\pi(i)} = u_i \circ \pi^{-1}$; 
		\item For each $i \in N$ and $\tau \in T_N$, $u_i = u_{\tau(i)} \circ \tau$; and
		\item For each $i \in N$ and $\tau \in T_N$, $u_i = u_{\tau(i)} \circ \tau^{-1}$.
	\end{enumerate}
	
	\begin{proof}		
		Condition (ii) follows trivially from Condition (i). Now suppose Condition (ii) is satisfied and let $H$ be a transitive subgroup of player permutations under which $\Gamma$ is invariant. Let $\pi \in S_N$, $i \in N$ and $\tau \in H$ such that $\tau(i) = \pi(i)$. Since $(\tau^{-1} \circ \pi) \in S_{N-\{i\}}$ it follows from weak anonymity that $u_i = u_i \circ (\tau^{-1} \circ \pi)$. It also follows from standard symmetry that $u_i = u_{\tau(i)} \circ \tau$, putting these two bits of information together we have $u_i = u_i \circ (\tau^{-1} \circ \pi) = (u_{\tau(i)} \circ \tau) \circ (\tau^{-1} \circ \pi) = u_{\tau(i)} \circ (\tau \circ \tau^{-1}) \circ \pi = u_{\tau(i)} \circ \pi = u_{\pi(i)} \circ \pi$.
	
		Suppose Condition (i) is satisfied, then for each $i \in N$ and $\pi \in S_N$, $u_{\pi(i)} = u_{\pi(i)} \circ (\pi \circ \pi^{-1}) = (u_{\pi(i)} \circ \pi) \circ \pi^{-1} = u_i \circ \pi^{-1}$. The converse works the same in reverse giving equivalence of Conditions (i) and (iii). 
		
		Condition (i) implies Condition (iv) since $T_N \subseteq S_N$, and Condition (iv) implies Condition (i) directly from Corollary \ref{utilityactionprop} and that $\langle{T_N}\rangle = S_N$. Conditions (iv) and (v) are equivalent since each transposition is its own inverse.
	\end{proof}
\end{theorem}

Condition (iii) in Theorem \ref{basicsymequivthm} was used by von Neumann and Morgenstern \cite{VNM}, which was ideal for their chosen notation of permutations acting on the right of players and strategy profiles. Of course any generating set of $S_N$ may replace $T_N$ in Condition (iv) of Theorem \ref{basicsymequivthm}.

It is worth noting that it is easy to mistakenly use the following inequivalent condition: for each $i \in N$ and $\pi \in S_N$, $u_i = u_{\pi(i)} \circ \pi^{-1}$ \cite[Definition 7]{DMaskin}. However this does not permute the players and strategy profiles correctly as the right hand side does not have player $\pi(i)$ playing the strategy that player $i$ is playing, which we illustrate using Example \ref{fullsymeg}. 

Let $\pi = (123) \in S_3$, the incorrect condition given in \cite[Definition 7]{DMaskin} requires that for each $i \in N$ and $(s_1, s_2, s_3) \in A$, we have $u_i(s_1, s_2, s_3) = u_{\pi(i)}(s_{\pi(1)}, s_{\pi(2)}, s_{\pi(3)}) = u_{\pi(i)}(s_2, s_3, s_1)$. By considering $(b, a, a) \in A$, we see that $3 = u_1(b, a, a) \neq u_2(a, a, b) = 2$. It should be fairly obvious that if we are mapping player 1 to player 2 and player 1 is playing $b$ then we want the mapped strategy profile to have player 2 playing $b$.

Since $T_N \subseteq S_N$, it follows from Condition (v) in Theorem \ref{basicsymequivthm} that the incorrect condition in \cite[Definition 7]{DMaskin} is somewhat surprisingly a more restrictive condition than the conditions in Theorem \ref{basicsymequivthm}. When $n=2$, since each transposition is its own inverse, the incorrect condition in \cite[Definition 7]{DMaskin} is equivalent to the conditions in Theorem \ref{basicsymequivthm}. We now establish that for $n \geq 3$ the incorrect condition in \cite[Definition 7]{DMaskin} is equivalent to the condition for a game being fully anonymous. 

\begin{lemma} \label{brandtlemma}
	\cite{brandt2009symmetries} The following conditions are equivalent:
	\begin{enumerate}
		\item $\Gamma$ is fully anonymous; and
		\item $\Gamma$ is fully symmetric and $u_i = u_j$ for all $i, j \in N$.
	\end{enumerate}
\end{lemma}

\begin{lemma} \label{DMlemma}
	Let $\pi, \tau \in S_N$. If $u_i = u_{\pi(i)} \circ \pi^{-1} = u_{\tau(i)} \circ \tau^{-1}$ for all $i \in N$ then $u_i = u_{(\tau \circ \pi)(i)} \circ (\pi \circ \tau)^{-1}$ for all $i \in N$.
	\begin{proof}
		For each $i \in N$, $u_i = u_{\pi(i)} \circ \pi^{-1} = (u_{\tau(\pi(i))} \circ \tau^{-1}) \circ \pi^{-1} = u_{(\tau \circ \pi)(i)} \circ (\pi \circ \tau)^{-1}$.
	\end{proof}
\end{lemma}

\begin{theorem} \label{DMprop}
	If $n \geq 3$ then the following conditions are equivalent:
	\begin{enumerate}
		\item $\Gamma$ is fully symmetric and $u_i = u_j$ for all $i, j \in N$; and
		\item For each $i \in N$ and $\pi \in S_N$, $u_i = u_{\pi(i)} \circ \pi^{-1}$.
	\end{enumerate}
	\begin{proof}
		Suppose Condition (i) holds, then for each $i \in N$ and $\pi \in S_N$, $u_i = u_{\pi^{-1}(i)} \circ \pi^{-1} = u_{\pi(i)} \circ \pi^{-1}$. Conversely suppose Condition (ii) holds, and hence that $\Gamma$ is fully symmetric. Let $i, j, k \in N$ be distinct. Since $(ik) \circ (ijk) \circ (jk) = (ijk)$ and $\bigl((jk) \circ (ijk) \circ (ik)\bigr)^{-1} = (ik) \circ (ikj) \circ (jk) =  e$, it follows from Lemma \ref{DMlemma} that $u_i = u_j$. 
	\end{proof}
\end{theorem}

We conclude this subsection by providing the reader with an accurate historical account of the mistake from \cite[Definition 7]{DMaskin} being identified. The mistake was first pointed out by the author with an edit on the 4\textsuperscript{th} of May 2011 to the Wikipedia page for symmetric games, which the author then revised on the 8\textsuperscript{th} of May 2011 due to not having a published reference for the author's claim that the definition is incorrect. Both of these edits are visible on the Wikipedia revision history for the symmetric games page \cite{WikiSGRV}. The mistake was also pointed out in the author's 2011 honours thesis \cite[Subsection 5.8]{ham2011honoursthesis}.

Upon contacting the authors from \cite{DMaskin} in 2018, the author received a response from Maskin suggesting that they made a slight mistake, unintentionally making the definition of symmetry given stronger than intended. Maskin suggested the mistake did not affect their own results, but has had the unfortunate effect of possibly leading other researchers astray. Prior to 2011 \cite{DMaskin} had 949 citations, and as at December 2018 it has 1,374 citations, so the author feels it is a good idea for the mistake to be noted to hopefully avoid any researchers being led astray in the future.

The mistake was also pointed out independently by Vester in his 2012 Masters thesis \cite[Appendix B]{vester2012symmetric}, who also proved the statement in Theorem \ref{DMprop}. Theorem \ref{DMprop} does not appear in the author's honours thesis, a proof was first released by the author publicly with the first revision of this paper uploaded to the arXiv November 2013, see \cite[Version 1]{ham2018arxivversion}. Hence credit goes to Vester for first releasing a proof publicly, see \cite[Theorem 32]{vester2012symmetric}.

Further, Tohm\'{e} et al. also proved the statement in Theorem \ref{DMprop} which they released in 2017, see \cite[Lemma 2.14]{tohme2019structural}.

\subsection{Notions of Fairness} \label{subsec:labeldepnotionsoffairness}
Interestingly, fairness has not appeared much in the game theory literature. Here we review where the term fair has appeared, introduce several new notions of fairness and begin examining how they relate to one another. Note that we will revisit these notions of fairness several times throughout the remainder of the paper.

\begin{definition}
	A $2$-player zero-sum game is \textit{fair} \cite[17.11, 28.1, 28.2]{VNM} if its value is $0$.
\end{definition}

\begin{proposition} 
	\cite{VNM} If a $2$-player zero-sum game is fully symmetric then it is fair. 
\end{proposition}

\begin{definition}
	We shall refer to a game $\Gamma = (N, A, u)$ as:
	\begin{enumerate}
		\item \textit{maximin fair} if $\underline{u}_i = \underline{u}_j$ for all $i, j \in N$;
		\item \textit{minimax fair} if $\overline{u}_i = \overline{u}_j$ for all $i, j \in N$;
		\item \textit{very-weakly-fair} if $\underline{u}_i = \underline{u}_j$ and $\overline{u}_i = \overline{u}_j$ for all $i, j \in N$;
		\item \textit{weakly-fair} if $\underline{u}_i = \overline{u}_i = \underline{u}_j = \overline{u}_j$ for all $i, j \in N$; 
		\item \textit{fair} if $\underline{u}_i = \overline{u}_i = 0$ for all $i \in N$;
		\item \textit{standard fair} if utility values are preserved under a transitive subgroup of the player permutations;
		\item \textit{fully fair} if utility values are preserved under all player permutations;
		\item \textit{standard ordinally fair} if there is a transitive subgroup of player permutations that preserve preferences over pure strategy profiles;
		\item \textit{fully ordinally fair} if all player permutations preserve preferences over pure strategy profiles;
		\item \textit{standard cardinally fair} if there is a transitive subgroup of player permutations that preserve preferences over mixed strategy profiles; and
		\item \textit{fully cardinally fair} if all player permutations preserve preferences over mixed strategy profiles.
	\end{enumerate}
\end{definition}

Note that ordinally symmetric games have been examined by Cao et al \cite{cao2018symmetric, cao2019ordinally}. Our definitions of standard and fully fair match up with our definitions of standard and fully symmetric games, a similar situation holds for the ordinal and cardinal definitions for ordinal and cardinal generalisations of our symmetric definitions to capture symmetry of payoff structure rather than symmetry of payoffs directly. 

An alternative way some people view fairness is requiring players reach equal payoffs under reasonable notions of perfect play. The author feels this would be more akin to the definition of symmetry in \cite{DMaskin}.

\begin{proposition} \label{prop:standsymgamesaremaxminandminmaxfair}
	If a (zero-sum) game $\Gamma = (N, A, u)$ is standard symmetric then it is maximin fair and minimax fair.
	
	\begin{proof}
		Since $\Gamma$ is standard symmetric there exists a transitive subgroup of the game invariants $T$ such that for each $\pi \in T$, $u_i = u_i \circ \pi$ for all $i \in N$. For each $i, j \in N$, since $T$ is transitive there exists $\pi \in T$ such that $\pi(i) = j$. Using $u_i = u_{\pi(i)} \circ \pi = u_j \circ \pi$, rearranging and then changing variables we have:
		\begin{align*}
			\text{(i) }\max_{\sigma_i \in \Delta(A_i)}\min_{\sigma_{-i} \in {\nabla(A)}_{-i}} u_i(\sigma_i, \sigma_{-i}) &= \max_{\sigma_i \in \Delta(A_i)}\min_{\sigma_{-i} \in \Delta(A)_{-i}} u_i(\sigma) \\
			&= \max_{\sigma_i \in \Delta(A_i)}\min_{\sigma_{-i} \in {\nabla(A)}_{-i}} (u_{\pi(i)} \circ \pi)(\sigma) \\
			&= \max_{\sigma_i \in \Delta(A_i)}\min_{\sigma_{-i} \in {\nabla(A)}_{-i}} u_j\left(\pi(\sigma)\right) \\
			&= \max_{\sigma_i \in \Delta(A_i)}\min_{\sigma_{-i} \in {\nabla(A)}_{-i}} u_j\left(\pi(\sigma)_j, \pi(\sigma)_{-j}\right) \\
			&= \max_{\sigma_i \in \Delta(A_i)}\min_{\sigma_{-i} \in {\nabla(A)}_{-i}} u_j\left(\pi(\sigma_i), \pi(\sigma_{-i})\right) \\
			&= \max_{\theta_j \in \Delta(A_j)}\min_{\theta_{-j} \in {\nabla(A)}_{-j}} u_j(\theta_j, \theta_{-j}); \text{ and} \\
			\text{(ii) }\min_{\sigma_{-i} \in {\nabla(A)}_{-i}}\max_{\sigma_i \in \Delta(A_i)} u_i(\sigma_i, \sigma_{-i}) &= \min_{\sigma_{-i} \in {\nabla(A)}_{-i}}\max_{\sigma_i \in \Delta(A_i)} u_i(\sigma) \\
			&= \min_{\sigma_{-i} \in {\nabla(A)}_{-i}}\max_{\sigma_i \in \Delta(A_i)} (u_{\pi(i)} \circ \pi)(\sigma) \\
			&= \min_{\sigma_{-i} \in {\nabla(A)}_{-i}}\max_{\sigma_i \in \Delta(A_i)} u_j\left(\pi(\sigma)\right) \\
			&= \min_{\sigma_{-i} \in {\nabla(A)}_{-i}}\max_{\sigma_i \in \Delta(A_i)} u_j\left(\pi(\sigma)_j, \pi(\sigma)_{-j}\right) \\
			&= \min_{\sigma_{-i} \in {\nabla(A)}_{-i}}\max_{\sigma_i \in \Delta(A_i)} u_j\left(\pi(\sigma_i), \pi(\sigma_{-i})\right) \\
			&= \min_{\theta_{-j} \in {\nabla(A)}_{-j}}\max_{\theta_j \in \Delta(A_j)} u_j(\theta_j, \theta_{-j}). 
		\end{align*} 
	\end{proof}
\end{proposition}

The reader may like to verify that Proposition \ref{prop:standsymgamesaremaxminandminmaxfair} holds in the label-independent case. It was established by von Neumann and Morgenstern \cite[Pages 165-166]{VNM} that every $2$-player standard symmetric zero-sum game is fair. The reader may also like to determine whether a zero-sum standard symmetric game $\Gamma = (N, A, u)$ is necessarily fair.

    \section{Morphisms Between Games} \label{sec:morphisms}
There are two important reasons why our simplifying assumption that players have the same strategy labels leaves our analysis incomplete. Our first reason is that relabelling the strategies for a standard symmetric game leads to a strategically equivalent game that may no longer be considered symmetric inside our label-dependent framework. 

Ideally we want to be able to determine when two games merely differ by player and strategy labels without having to go through and check all possible rearrangements of the labels.

Our second reason is that there are weaker notions of fairness that cannot be captured within our label-dependent framework. As a motivating example consider Matching Pennies.

\begin{example} \label{MPeg}
	Matching Pennies
	\begin{center}
		\begin{game}{2}{2}
			      \> $H$    \> $T$ \\
			$H$   \> $1,-1$  \> $-1,1$ \\
			$T$   \> $-1,1$  \> $1,-1$
		\end{game} 
	\end{center}
	It is clear just by looking at the payoff matrix that Matching Pennies is fair, yet inside our label dependent framework the only invariant is the identity permutation, a problem that persists if we swap the strategy labels for either or both of the players.             
\end{example}

\subsection{Game Bijections}    
	\begin{definition}
		A \textit{game bijection} from $\Gamma_1 = (N, A, u)$ to $\Gamma_2 = (M, B, v)$ consists of a bijection $\pi:N\rightarrow M$ and for each player $i \in N$, a bijection $\tau_i:A_i\rightarrow B_{\pi(i)}$, which we denote as $\bigl(\pi; (\tau_i)_{i \in N}\bigr)$.
	\end{definition}

	More on game bijections can be found in \cite{IsoComplexity}. We denote the set of game bijections from $\Gamma_1$ to $\Gamma_2$ as $\bij(\Gamma_1, \Gamma_2)$, or simply $S_{\Gamma}$ for the bijections from a game $\Gamma$ to itself. Let $g = \bigl(\pi; (\tau_i)_{i \in N}\bigr) \in \bij(\Gamma_1, \Gamma_2)$, $i \in N$, $s_i \in A_i$ and $s \in A$, using similar notation to our label-dependent framework we denote $\pi(i)$ as $g(i)$, $\tau_i(s_i)$ as $g(s_i)$, $\bigl(\tau_{\pi^{-1}(j)}(s_{\pi^{-1}(j)})\bigr)_{j \in M} \in B$ as $g(s)$ giving $\bigl(g(s)\bigr)_{g(i)} = \tau_i(s_i) = g(s_i)$, and the map $s \mapsto u_{g(i)}\bigl(g(s)\bigr)$ as $u_{g(i)} \circ g$. 	
	
	\begin{example} \label{egisomgames}
		Consider the following $2$-player games.            
        \begin{center}
            \begin{game}{2}{2}[$\Gamma_1$]
                        \> $c$  \> $d$ \\
                $a$   \> $1,2$  \> $3,4$ \\
                $b$   \> $5,6$  \> $7,8$
            \end{game}
            \hspace*{10mm} 
            \begin{game}{2}{2}[$\Gamma_2$]
                        \> $h$ \> $i$ \\
                $e$   \> $4,3$ \> $8,7$ \\
                $f$   \> $2,1$ \> $6,5$ \\
            \end{game} 
        \end{center}
        
		Given $(a, c) \in A$ and $g = \bigl((12) ; \bigl(\begin{smallmatrix} a & b \\ h & i \end{smallmatrix}\bigr), \bigl(\begin{smallmatrix} c & d \\ f & e \end{smallmatrix}\bigr)\bigr) \in \bij(\Gamma_1, \Gamma_2)$, $g(a,c) = (f,h)$.
	\end{example}
        
	Let $\Gamma_3 = (L, C, w)$ also be a game. For $g = \bigl(\pi; (\tau_i)_{i \in N}\bigr) \in \bij(\Gamma_1, \Gamma_2)$ and $h = \bigl(\eta; (\phi_j)_{j \in M}\bigr) \in \bij(\Gamma_2, \Gamma_3)$, their \textit{composite}, denoted $h\circ g$, is $\bigl(\eta\circ\pi; (\phi_{\pi(i)}\circ\tau_i)_{i \in N}\bigr) \in \bij(\Gamma_1, \Gamma_3)$, and the \textit{inverse} of $g$, denoted $g^{-1}$, is $\bigl(\pi^{-1}; (\tau^{-1}_{\pi^{-1}(j)})_{j \in M}\bigr) \in \bij(\Gamma_2, \Gamma_1)$.   
   	
   	\begin{example}
   		Consider Example \ref{stdsymeg} except with strategy labels $A_2 =\{c, d\}$ and $A_3 = \{e, f\}$ for players $2$ and $3$ respectively. We compose and invert bijections $g = \bigl((123) ; \bigl(\begin{smallmatrix} a & b \\ d & c \end{smallmatrix}\bigr), \bigl(\begin{smallmatrix} c & d \\ e & f \end{smallmatrix}\bigr), \bigl(\begin{smallmatrix} e & f \\ b & a \end{smallmatrix}\bigr)\bigr)$, $h = \bigl((12) ; \bigl(\begin{smallmatrix} a & b \\ c & d \end{smallmatrix}\bigr), \bigl(\begin{smallmatrix} c & d \\ a & b \end{smallmatrix}\bigr), \bigl(\begin{smallmatrix} e & f \\ f & e \end{smallmatrix}\bigr)\bigr) \in S_{\Gamma}$ as follows:
   		\begin{align*}
   			h \circ g &= \bigl((12) ; \bigl(\begin{smallmatrix} a & b \\ c & d \end{smallmatrix}\bigr), \bigl(\begin{smallmatrix} c & d \\ a & b \end{smallmatrix}\bigr), \bigl(\begin{smallmatrix} e & f \\ f & e \end{smallmatrix}\bigr)\bigr) \circ \bigl((123) ; \bigl(\begin{smallmatrix} a & b \\ d & c \end{smallmatrix}\bigr), \bigl(\begin{smallmatrix} c & d \\ e & f \end{smallmatrix}\bigr), \bigl(\begin{smallmatrix} e & f \\ b & a \end{smallmatrix}\bigr)\bigr) \\
   				&= \bigl((12) \circ (123) ; \bigl(\begin{smallmatrix} c & d \\ a & b \end{smallmatrix}\bigr) \circ \bigl(\begin{smallmatrix} a & b \\ d & c \end{smallmatrix}\bigr), \bigl(\begin{smallmatrix} e & f \\ f & e \end{smallmatrix}\bigr) \circ \bigl(\begin{smallmatrix} c & d \\ e & f \end{smallmatrix}\bigr), \bigl(\begin{smallmatrix} a & b \\ c & d \end{smallmatrix}\bigr) \circ \bigl(\begin{smallmatrix} e & f \\ b & a \end{smallmatrix}\bigr)\bigr)\\
   				&= \bigl((23) ; \bigl(\begin{smallmatrix} a & b \\ b & a \end{smallmatrix}\bigr), \bigl(\begin{smallmatrix} c & d \\ f & e \end{smallmatrix}\bigr), \bigl(\begin{smallmatrix} e & f \\ d & c \end{smallmatrix}\bigr)\bigr) \text{; and} \\
   			g^{-1} &= \bigl((123) ; \bigl(\begin{smallmatrix} a & b \\ d & c \end{smallmatrix}\bigr), \bigl(\begin{smallmatrix} c & d \\ e & f \end{smallmatrix}\bigr), \bigl(\begin{smallmatrix} e & f \\ b & a \end{smallmatrix}\bigr)\bigr)^{-1} \\
   				   &= \bigl((123)^{-1} ; \bigl(\begin{smallmatrix} e & f \\ b & a \end{smallmatrix}\bigr)^{-1}, \bigl(\begin{smallmatrix} a & b \\ d & c \end{smallmatrix}\bigr)^{-1}, \bigl(\begin{smallmatrix} c & d \\ e & f \end{smallmatrix}\bigr)^{-1}\bigr)\\
   				   &= \bigl((132) ; \bigl(\begin{smallmatrix} a & b \\ f & e \end{smallmatrix}\bigr), \bigl(\begin{smallmatrix} c & d \\ b & a \end{smallmatrix}\bigr), \bigl(\begin{smallmatrix} e & f \\ c & d \end{smallmatrix}\bigr)\bigr).
   		\end{align*}
   	\end{example}
        
	\begin{lemma} 
		$(h \circ g)(s) = h(g(s))$ for all $s \in A$.
		\begin{proof}
			\begin{align*}
				(h \circ g)(s) &= \bigl(\eta \circ \pi; (\phi_{\pi(i)} \circ \tau_i)_{i \in N}\bigr)(s) \\
				&= \bigl(\phi_{\eta^{-1}(k)} \circ \tau_{(\eta \circ \pi)^{-1}(k)}(s_{(\eta \circ \pi)^{-1}(k)})\bigr)_{k \in L} \\
				&= \Bigl(\phi_{\eta^{-1}(k)}\bigl(\tau_{\pi^{-1}(\eta^{-1}(k))}(s_{\pi^{-1}(\eta^{-1}(k))})\bigr)\Bigr)_{k \in L} \\
				&= \Bigl(\phi_{\eta^{-1}(k)}\bigl(g(s)_{\eta^{-1}(k)}\bigr)\Bigr)_{k \in L} \\
				&= \Bigl(h\bigl(g(s)\bigr)_k\Bigr)_{k \in L} \\
				&= h(g(s)).
			\end{align*}
		\end{proof}
	\end{lemma}
	
	\begin{corollary}
		$u_{(h \circ g)(i)} \circ (h \circ g) = (u_{h(g(i))} \circ h) \circ g$ for all $i \in N$.
		
		\begin{proof}
			This follows identically to the proof of Corollary \ref{utilityactionprop}.
		\end{proof}
	\end{corollary}
        
    \begin{theorem} \label{bijgroupoidthm}
        Game bijections form a groupoid.
            
        \begin{proof}
            Let $\Gamma_3 = (P, C)$, $\Gamma_4 = (Q, D)$, $f = \bigl(\pi; (\tau_i)_{i \in N}\bigr) \in \bij(\Gamma_1, \Gamma_2)$, $g = \bigl(\eta; (\phi_j)_{j \in M}\bigr) \in \bij(\Gamma_2, \Gamma_3)$, $h = \bigl(\xi ; (\lambda_k)_{k \in P}\bigr) \in \bij(\Gamma_3, \Gamma_4)$. Then:
            \begin{align*}
                f \circ \text{id}_{\Gamma_1} &= \bigl(\pi \circ \text{id}_N; (\tau_i \circ \text{id}_{A_i})_{i \in N}\bigr) \\ 
                      &= f = \bigl(\text{id}_M \circ \pi; (\text{id}_{B_{\pi(i)}} \circ \tau_i)_{i \in N}\bigr) = \text{id}_{\Gamma_2} \circ f; \\
                f \circ f^{-1} &= \bigl(\pi \circ \pi^{-1}; (\tau_{\pi^{-1}(j)} \circ \tau^{-1}_{\pi^{-1}(j)})_{j \in M}\bigr) = \text{id}_{\Gamma_2}; \\
                f^{-1} \circ f &= \bigl(\pi^{-1} \circ \pi; (\tau^{-1}_{\pi^{-1}(\pi(i))} \circ \tau_i)_{i \in N}\bigr) = \text{id}_{\Gamma_1}; \text{ and} \\
                h \circ (g \circ f) &= \bigl(\xi ; (\lambda_k)_{k \in P}\bigr) \circ \bigl(\eta \circ \pi; (\phi_{\pi(i)} \circ \tau_i)_{i \in N}\bigr) \\
                    &= \bigl(\xi \circ \eta \circ \pi; (\lambda_{(\eta \circ \pi)(i)} \circ \phi_{\pi(i)} \circ \tau_i)_{i \in N}\bigr) \\
                    &= \bigl(\xi \circ \eta; (\lambda_{\eta(j)} \circ \phi_j)_{j \in M}\bigr) \circ \bigl(\pi; (\tau_i)_{i \in N}\bigr) = (h \circ g) \circ f.
			\end{align*}
		\end{proof}
	\end{theorem}

    \subsection {Game Isomorphisms}
	Game isomorphisms are game bijections that preserve strategic structure, they are useful for establishing strategic equivalence between games, or as we will be using them, for considering label-independent notions of symmetry.
        
	We will only require the strictest notion of game isomorphism to explore label-independent notions of symmetry, treating two games as isomorphic when they differ only by the player and strategy labels. However one can define ordinal and cardinal game isomorphisms by requiring preservation of preferences over pure and mixed strategy profiles respectively, then characterise each by the existence of increasing monotonic and affine transformations respectively, see \cite[Propositions 4.3.2 and 4.3.5]{ham2011honoursthesis}. A discussion on the computational complexity of deciding whether two games satisfy various notions of equivalence can be found in Gabarr\'{o} et al. \cite{IsoComplexity}.
        
	\begin{definition}
		A bijection $g \in \bij(\Gamma_1, \Gamma_2)$ is a \textit{game isomorphism} if $u_i = v_{g(i)}\circ g$ for all $i \in N$. 
	\end{definition}

	We denote by $\isom(\Gamma_1, \Gamma_2)$ the set of isomorphisms from $\Gamma_1$ to $\Gamma_2$, and write $\Gamma_1 \cong \Gamma_2$ when $\isom(\Gamma_1, \Gamma_2)$ is non-empty. The reader may like to verify that the bijection in Example \ref{egisomgames} is in fact an isomorphism. For example, $u_1(a, d) = v_{g(1)}\bigl(g(a, d)\bigr) = v_2(e, h)$.
            
	\begin{theorem} 
		Game isomorphisms form a groupoid.
            
		\begin{proof}                
			For each $g \in \isom(\Gamma_1, \Gamma_2)$ and $j \in M$, $v_j = (v_j \circ g) \circ g^{-1} = u_{g^{-1}(j)} \circ g^{-1}$, giving us $g^{-1} \in \isom(\Gamma_2, \Gamma_1)$. Let $\Gamma_3 = (P, C, w)$, then for each $g \in \isom(\Gamma_1, \Gamma_2)$, $h \in \isom(\Gamma_2, \Gamma_3)$ and $i \in N$, $u_i = v_{g(i)} \circ g = (w_{h(g(i))} \circ h) \circ g = w_{(h \circ g)(i)} \circ (h \circ g)$, giving us $(h \circ g) \in \isom(\Gamma_1, \Gamma_3)$.
                  
			The remaining conditions follow from Theorem \ref{bijgroupoidthm}.
		\end{proof}
	\end{theorem}
        
	\begin{corollary} \label{isocorollary}
		If $\Gamma_1 \cong \Gamma_2 \cong \Gamma_3$ then $\isom(\Gamma_1, \Gamma_2) \cong \isom(\Gamma_2, \Gamma_3)$.
	\end{corollary}
        
	Game isomorphisms induce an equivalence relation where games in the same equivalence class have the same strategic structure. There is a finite number of ordinal equivalence classes for games with both a fixed number of players and fixed number of strategies for each of the players. Goforth and Robinson \cite{GoforthRobinson} counted 144 ordinal equivalence classes for the 2-player 2-strategy games.
        
\subsection{Bijections Acting on Strategy Profiles}
The bijections $S_{\Gamma}$ from a game to itself form a group that acts on the players and strategy profiles. In fact for an $m$-strategy game $S_{\Gamma}$ is isomorphic to the wreath product $S_N \wr S_M$ where $M = \{1, \ldots, m\}$, which may be seen by setting $A_i = M$ for all $i \in N$. 

Given a game bijection $g = \bigl(\pi; (\tau_i)_{i \in N}\bigr) \in S_{\Gamma}$, we refer to $\pi$ as \textit{the player permutation used by $g$} and say that two game bijections $g, h \in S_{\Gamma}$ \textit{have the same player permutation} if the player permutations used by $g$ and $h$ are identical. 

Let $G$ be a subgroup of $S_{\Gamma}$. We denote the subgroup of player permutations used by game bijections in $G$ as $\overrightarrow{G}$. Furthermore, we say that $G$ is \textit{player transitive} if $G$ acts transitively on $N$, \textit{player $n$-transitive} if $G$ acts $n$-transitively on $N$, and \textit{only-transitive} if $G$ acts transitively and not $n$-transitively on $N$.

\begin{lemma} \label{cosetprop}
	Two bijections $g, h \in G$ have the same player permutation if and only if they are in the same coset of $G/G_N$.
	\begin{proof}
		Suppose $g , h$ have the same player permutation, then $h = g \circ (g^{-1} \circ h) \in (g \circ G_N)$. The converse is obvious.
	\end{proof}
\end{lemma}

Hence the factor group $G/G_N$ merely tells us what player permutations are used by $G$. 

\begin{corollary} \label{cosetcor}
	$G/G_N \cong \overrightarrow{G}$.
\end{corollary}

The isomorphisms from a game to itself form a subgroup of the game bijections called the \textit{automorphism group} of $\Gamma$, which we denote as $\Aut(\Gamma)$. Game automorphisms capture the notion of players being indifferent between the current positions and an alternative arrangement of positions. Note our definition is equivalent to the definition used by Nash \cite{NashNCG}.

For the sake of brevity, we refer to a subgroup of $\Aut(\Gamma)$ as a subgroup of $\Gamma$, denote the stabiliser subgroup of $\Aut(\Gamma)$ on $N$ as $\Gamma_N$, and denote the player permutations used by $\Aut(\Gamma)$ as $\overrightarrow{\Gamma}$.

\subsection{Strategy Triviality and Matchings}
Now that players need not have the same strategy labels, we seek a way to determine which subgroups of $S_{\Gamma}$ act on strategy profiles in an equivalent way to permutations for some relabelling of the strategies. Stein \cite{NoahXE} introduced strategy triviality for this purpose.

\begin{definition}
	A subgroup $G$ of $S_{\Gamma}$ is \textit{strategy trivial} \cite{NoahXE} if for each $i \in N$, $g(s_i) = s_i$ for all $g \in G_i$ and $s_i \in A_i$. 
\end{definition}

\begin{lemma} \label{NoahLemma} \cite{NoahXE}
	If $G$ is strategy trivial then for each $g, h \in G$ such that $g(i) = h(i)$, $g(s_i) = h(s_i)$ for all $s_i \in A_i$.
	\begin{proof}
		Since $(g^{-1} \circ h) \in G_i$, by strategy triviality, $g(s_i) = g\bigl((g^{-1} \circ h)(s_i)\bigr) = (g \circ g^{-1})\bigl(h(s_i)\bigr) = h(s_i)$.
	\end{proof}
\end{lemma}

\begin{corollary} \label{cosetprop2}
	If $G$ is strategy trivial then $G_N = \{\id_{\Gamma}\}$.
\end{corollary}

Hence strategy trivial subgroups have at most one bijection for each player permutation. Example \ref{egntransstandnonfull} establishes that the converse of Corollary \ref{cosetprop2} is false.

\begin{corollary} \label{NoahCor} 
	If $G \leq S_{\Gamma}$ is strategy trivial then for each $i \in N$ and $\tau \in \overrightarrow{G}$, there exists $g_{i\tau(i)} \in \bij(A_i, A_{\tau(i)})$ such that $G = \{(\pi; (g_{i\pi(i)})_{i \in N}): \pi \in \overrightarrow{G}\}$.
\end{corollary}

It follows that all paths from one player to another map the strategies in a canonical manner. Hence if $G$ is also player transitive then the strategy sets are matched such that they can be treated as the same set. We now introduce \textit{matchings} to formalise what is meant by the strategy sets being matched.

\begin{definition}
	A \textit{matching of $A_1, \ldots, A_n$} is a relation $M \subseteq \times_{i \in N} A_i$ which is $i$-total and $i$-unique for all $i \in N$. 
\end{definition}

\begin{example} \label{matchingeg}
Let $A_1 = \{a, b\}$, $A_2 = \{c, d\}$ and $A_3 = \{e, f\}$. One matching of $A_1 \times A_2 \times A_3$ is $M = \{(a,d,f), (b,c,e)\}$.
\begin{center}
\begin{tikzpicture}
	nodes={draw, ultra thick}
	\node (a) at (0,0) {$a$};	\node (c) at (1,0) {$c$};	\node (e) at (2,0) {$e$};
	\node (b) at (0,-0.5) {$b$};	\node (d) at (1,-0.5) {$d$};	\node (f) at (2,-0.5) {$f$};
	\draw (a) -- (d) -- (f);
	\draw (b) -- (c) -- (e);
\end{tikzpicture}
\end{center}
\end{example}

From a game theoretic point of view, a matching is a subset $M$ of the strategy profiles where for each $i \in N$ and $a_i \in A_i$ there is exactly one $s \in M$ such that $s_i = a_i$, and hence $|M| = m$. 

For each $i, j \in N$, a matching $M$ induces a bijection $M_{ij} \in \bij(A_i, A_j)$ where, given $a_i \in A_i$, $M_{ij}(a_i)$ is the unique $a_j \in A_j$ such that there exists $s \in M$ with $s_i = a_i$ and $s_j = a_j$. For example given the matching in Example \ref{matchingeg}, $M_{31} = \bigl(\begin{smallmatrix} e & f \\ b & a \end{smallmatrix}\bigr)$.

\begin{lemma}
	$\{M_{ij}: i, j \in N\}$ is a groupoid. 
	\begin{proof}
		It follows by definition that for each $i, j, k \in N$, $M_{ii} = \id_{A_i}$, $M_{ij}^{-1} = M_{ji}$ and $M_{jk} \circ M_{ij} = M_{ik}$. Now for each $i, j, k, l \in N$, $M_{ij} \circ M_{ii} = M_{ij} = M_{jj} \circ M_{ij}$, $M_{kl} \circ (M_{jk} \circ M_{ij}) = M_{kl} \circ M_{ik} = M_{il} = M_{jl} \circ M_{ij} = (M_{kl} \circ M_{jk}) \circ M_{ij}$, $M_{ij} \circ M_{ij}^{-1} = M_{ij} \circ M_{ji} = M_{jj}$ and $M_{ij}^{-1} \circ M_{ij} = M_{ji} \circ M_{ij} = M_{ii}$.
	\end{proof}
\end{lemma}

Furthermore, for each $\pi \in S_N$, a matching $M$ induces a game bijection $\bigl(\pi; (M_{i\pi(i)})_{i \in N}\bigr) \in S_{\Gamma}$, which we denote as $M_{\pi}$. For example given the matching in Example \ref{matchingeg}, $M_{(13)} = \bigl((13) ; \bigl(\begin{smallmatrix} a & b \\ f & e \end{smallmatrix}\bigr), \bigl(\begin{smallmatrix} c & d \\ c & d \end{smallmatrix}\bigr), \bigl(\begin{smallmatrix} e & f \\ b & a \end{smallmatrix}\bigr)\bigr)$. 

For each $H \subseteq S_N$ and matching $M$ we denote the set $\{M_{\pi}: \pi \in H\}$ of bijections induced by $H$ as $M_H$. For example given a subgroup $G$ of $S_{\Gamma}$ we have $M_{\overrightarrow{G}} = \{M_{\pi}: \pi \in \overrightarrow{G}\}$. 

\begin{lemma} \label{matchinghomoprop}
	$M:S_N\rightarrow{S_{\Gamma}}$ is a homomorphism.
	\begin{proof}
		Let $\pi, \phi \in S_N$, then $M_{\phi} \circ M_{\pi} = \bigl(\phi; (M_{i\phi(i)})_{i \in N}\bigr) \circ \bigl(\pi; (M_{i\pi(i)})_{i \in N}\bigr) = \bigl(\phi \circ \pi; (M_{\pi(i)(\phi\circ\pi)(i)} \circ M_{i\pi(i)})_{i \in N}\bigr) \newline = \bigl(\phi \circ \pi; (M_{i(\phi\circ\pi)(i)})_{i \in N}\bigr) = M_{(\phi \circ \pi)}$.
	\end{proof}
\end{lemma}

\begin{corollary}
	$M_{\pi^{-1}} = M_{\pi}^{-1}$ for all $\pi \in S_N$.
\end{corollary}

\begin{lemma} \label{matchingfixedpointlemma}
	For each $\pi \in S_N$, $M_{\pi}(s) = s$ for all $s \in M$.
	\begin{proof}
		For each $i \in N$, $\bigl(M_{\pi}(s)\bigr)_i = M_{\pi^{-1}(i)i}(s_{\pi^{-1}(i)}) = s_i$.
	\end{proof}
\end{lemma}

If we relabel the strategies played in each $s \in M$ to be the same, giving players the same strategy labels, then each permutation $\pi \in S_N$ acts on our relabelled strategy profiles equivalently to how $M_{\pi}$ acts on our original strategy profiles. Hence a subgroup $G$ of $S_{\Gamma}$ acts on strategy profiles equivalently to permutations for some relabelling of the strategies precisely when $G = M_{\overrightarrow{G}}$ for some matching $M$, which we now establish occurs precisely when $G$ is strategy trivial. 

\begin{theorem} \label{strattrivmatchingthm}
	Let $G \leq S_{\Gamma}$ be player transitive. There exists a matching $M$ such that $M_{\overrightarrow{G}} = G$ if and only if $G$ is strategy trivial.
	\begin{proof}
		Suppose there exists a matching $M$ such that $M_{\overrightarrow{G}} = G$. That $M_{\overrightarrow{G}} \leq S_{\Gamma}$ follows from Lemma \ref{matchinghomoprop}. Now for each $i \in N$ and $g \in G_i$, $M_{ig(i)} = M_{ii} = \id_{A_i}$.
		
		Conversely suppose $G$ is strategy trivial. By Corollary \ref{NoahCor}, for each $i \in N$ and $\tau \in \overrightarrow{G}$ there exists $g_{i\tau(i)} \in \bij(A_i, A_{\tau(i)})$ such that $G = \{\bigl(\pi; (g_{i\pi(i)})_{i \in N}\bigr): \pi \in \overrightarrow{G}\}$. 
		
		Let $i \in N$ and $M = \{(g_{ij}(a_i))_{j \in N}: a_i \in A_i\}$. $M$ is a matching since for each $j \in N$ and $a_j \in A_j$, there exists a unique strategy $a_i \in A_i$ for player $i$ such that $g_{ij}(a_i) = a_j$. Furthermore $M$ is independent of $i$ since for each $k \in N$, $\bigl(g_{ij}(a_i)\bigr)_{j \in N} = \bigl((g_{kj} \circ g_{ik})(a_i)\bigr)_{j \in N}$. Hence $M_{kl} = g_{kl}$ for all $k, l \in N$, giving us $M_{\pi} = \bigl(\pi; (M_{i\pi(i)})_{i \in N}\bigr) = \bigl(\pi; (g_{i\pi(i)})_{i \in N}\bigr) \in G$ for all $\pi \in \overrightarrow{G}$.
	\end{proof} 
\end{theorem}

Hence weakly anonymous games may be characterised as follows, similarly for anonymous and fully anonymous games.

\begin{corollary}
	The following conditions are equivalent:
	\begin{enumerate}
		\item There exists weakly anonymous $\Gamma'$ such that $\Gamma \cong \Gamma'$;
		\item There exists player $n$-transitive and strategy trivial $G \leq \Gamma$ such that for each $i \in N$ and $g \in G_i$, $u_i = u_i \circ g$; and
		\item There exists a matching $M$ such that for each $i \in N$ and $\pi \in S_{N-\{i\}}$, $u_i = u_i \circ M_{\pi}$.
	\end{enumerate}
\end{corollary}

We denote by $M(n, m)$ the set of matchings for an $n$-player $m$-strategy game. 

\begin{example}
	\begin{enumerate}
		\item If $m = n = 2$ then, letting $A_1 = \{a, b\}$ and $A_2 = \{c, d\}$, 
		\begin{align*}
			M(2, 2) = \bigl\{\{(a, c), (b, d)\}, \{(a, d), (b, c)\}\bigr\}.
		\end{align*}
		\item If $m = 3$ and $n = 2$ then, letting $A_1 = \{a, b, c\}$ and $A_2 = \{d, e, f\}$, \
		\begin{align*}
			M(2, 3) = \bigl\{&\{(a, d), (b, e), (c, f)\}, \{(a, d), (b, f), (c, e)\}, \{(a, e), (b, d), (c, f)\}, \\
						&\{(a, e), (b, f), (c, d)\}, \{(a, f), (b, d), (c, e)\}, \{(a, f), (b, e), (c, d)\}\bigr\}.
		\end{align*}
	\end{enumerate}
\end{example}

There are a number of ways to count the number of matchings in $M(n, m)$. Below we present one, though note an alternative is to establish that $M(n, m) \cong \bij(A_1, A_2) \times \ldots \times \bij(A_{n-1}, A_n)$. 

\begin{lemma} \label{matchinglemma1}
	For each $n \geq 2$: $M(n, 2)$ is a partition of $A$; and $|M(n, 2)| = 2^{n-1}$.
	\begin{proof}
		For each $s \in A$, the profile $s'$ where each player swaps their strategy choice is the unique profile in $A$ such that $\{s, s'\} \in M(n, 2)$. Consequently $|M(n, 2)| = \frac{|A|}{2} = 2^{n-1}$. 
	\end{proof}
\end{lemma}

\begin{lemma} \label{matchinglemma2}
	For each $n \geq 2$ and $m \geq 3$, $|M(n, m)| = m^{n-1}|M(n, m-1)|$.
	\begin{proof}
		Let $i \in N$. Each $a_i$ can be matched with each $a_{-i} \in A_{-i}$ and $ |A_{-i}| = m^{n-1}$. Furthermore, for each $(a_i, a_{-i})$ there are $|M(n, m-1)|$ ways to match the remaining $m-1$ strategies of the $n$ players.
	\end{proof}
\end{lemma}

\begin{theorem}
	For each $m, n \geq 2$, $|M(n, m)| = (m!)^{n-1}$.
	\begin{proof}
		This follows inductively from Lemmas \ref{matchinglemma1} and \ref{matchinglemma2}.
	\end{proof}
\end{theorem}

    \section{Label-Independent Notions of Symmetry} \label{sec:labindnotions}
\subsection{Notions of Symmetry} \label{subsec:labelindepnotionsofsymmetry}
Similar to our label-independent characterisations of our label-dependent notions of anonymity, Theorem \ref{strattrivmatchingthm} gives us the following label-independent characterisations of our label-dependent notions of fairness.

\begin{corollary} \label{indchar1}
	The following conditions are equivalent:
	\begin{enumerate}
		\item There exists standard symmetric $\Gamma'$ such that $\Gamma \cong \Gamma'$;
		\item $\Gamma$ has a player transitive and strategy trivial subgroup $G$; and
		\item There exists a matching $M$ and player transitive $T \leq S_N$ such that $M_{T} \leq \Aut(\Gamma)$.
	\end{enumerate}
\end{corollary}

\begin{corollary}  \label{indchar2}
	The following conditions are equivalent:
	\begin{enumerate}
		\item There exists fully symmetric $\Gamma'$ such that $\Gamma \cong \Gamma'$;
		\item $\Gamma$ has a player $n$-transitive and strategy trivial subgroup $G$; and
		\item There exists a matching $M$ such that $M_{S_N} \leq \Aut(\Gamma)$.
	\end{enumerate}
\end{corollary}

Henceforth we will use fully and standard symmetric to refer to our label-independent characterisations. 

\begin{corollary}
	If $\Gamma$ is standard symmetric then there exists a matching $M$ such that for each $s \in M$, $u_i(s) = u_j(s)$ for all $i, j \in N$.
	\begin{proof}
		This follows from Lemma \ref{matchingfixedpointlemma}.
	\end{proof}
\end{corollary}

Remember that the defining features for standard and fully symmetric games inside our label-dependent framework were that players be indifferent between which position they play and the arrangement of the players respectively. Inside our label-independent framework, these defining features capture larger classes of fair games.

\begin{definition}
	A game is \textit{symmetric} \cite{NoahXE} if its automorphism group is player transitive and \textit{$n$-transitively symmetric} if its automorphism group is player $n$-transitive.
\end{definition}

\begin{example} The automorphism group of Matching Pennies in Example \ref{MPeg} is \begin{align*}
		\Aut(\Gamma) = \langle &\bigl((12) ; \bigl(\begin{smallmatrix} H & T \\ H & T \end{smallmatrix}\bigr), \bigl(\begin{smallmatrix} H & T \\ T & H \end{smallmatrix}\bigr)\bigr)\rangle \\
		= \{ &\bigl(e ; \bigl(\begin{smallmatrix} H & T \\ H & T \end{smallmatrix}\bigr), \bigl(\begin{smallmatrix} H & T \\ H & T \end{smallmatrix}\bigr)\bigr), 
		\bigl(e ; \bigl(\begin{smallmatrix} H & T \\ T & H \end{smallmatrix}\bigr), \bigl(\begin{smallmatrix} H & T \\ T & H \end{smallmatrix}\bigr)\bigr), \\
		&\bigl((12) ; \bigl(\begin{smallmatrix} H & T \\ H & T \end{smallmatrix}\bigr), \bigl(\begin{smallmatrix} H & T \\ T & H \end{smallmatrix}\bigr)\bigr),
		\bigl((12) ; \bigl(\begin{smallmatrix} H & T \\ T & H \end{smallmatrix}\bigr), \bigl(\begin{smallmatrix} H & T \\ H & T \end{smallmatrix}\bigr)\bigr) \}.
	\end{align*}
	Since $\Aut(\Gamma)$ is player transitive, is not strategy trivial and contains no proper transitive subgroups, Matching Pennies is an $n$-transitively non-standard symmetric game.
\end{example}

Peleg et al. \cite{peleg1999canonical} and Sudh\"{o}lter et al. \cite{sudholter2000canonical} define a game to be symmetric if $\Aut(\Gamma)/\Gamma_N \cong S_N$. It follows immediately from Corollary \ref{cosetcor} that this is equivalent to a game being $n$-transitively symmetric, and furthermore that $\Aut(\Gamma)/\Gamma_N$ being isomorphic to some transitive subgroup of $S_N$ is equivalent to a game being symmetric.

We now consider games which have a subgroup $G$ isomorphic to $S_N$ with $G_N = \{\id_{\Gamma}\}$. Fully symmetric games obviously satisfy this condition, Example \ref{egntransstandnonfull} shows that the converse of this is false. Below we show that all games satisfying this condition are $n$-transitively standard symmetric games; the author has been unable to show whether the converse holds.

\begin{proposition} \label{subisoimpliesinter}
	If $\Gamma$ has a subgroup $G$ isomorphic to $S_N$ with $G_N = \{\id_{\Gamma}\}$ then it is $n$-transitively standard symmetric. 
	\begin{proof}
		$n$-transitivity of $\Gamma$ follows from $\overrightarrow{G} = S_N$. Now since each $n$-cycle generates a regular subgroup of $S_N$, the subgroup of $G$ generated by an automorphism whose player permutation is an $n$-cycle is transitive and strategy trivial, hence $\Gamma$ is standard symmetric.
	\end{proof}
\end{proposition}

We end our exploration of symmetry notions with games that have a transitive subgroup $G$ isomorphic to $\overrightarrow{\Gamma}$ with $G_N = \{\id_{\Gamma}\}$. Standard symmetric games obviously satisfy this condition. To look at the converse we consider the argument used in Proposition \ref{subisoimpliesinter}. 

If all transitive subgroups of $S_N$ had regular subgroups then games with a transitive subgroup $G$ isomorphic to $\overrightarrow{\Gamma}$ with $G_N = \{\id_{\Gamma}\}$ would be standard symmetric. However this is not the case, Hulpke \cite{hulpke2005constructing} listed the non-regular minimally transitive permutation subgroups up to degree $30$. The smallest example is $\langle (14) \circ (25), (135) \circ (246)\rangle$ of degree $6$ and order $12$. For more information on how the transitive subgroups are constructed see for example \cite[Algorithm 8.1]{hulpke2016connected}. There is a GAP \cite{GAP4} library for transitive groups by Hulpke with manual \cite{hulpketransgrp}.

We will see in Example \ref{sixplayereg} that games which have a transitive subgroup $G$ isomorphic to $\overrightarrow{\Gamma}$ with $G_N = \{\id_{\Gamma}\}$ need not be standard symmetric. 

\subsection{Fairness Discussion} \label{subsec:fairnessdiscussion}
So far we have considered many notions of symmetry, see Subsections \ref{subsec:labeldepnotionsofsymmetry} and \ref{subsec:labelindepnotionsofsymmetry}, which we have also been referring to as notions of fairness. 

When the author was explaining the similarities between symmetry and fairness in the context of games to James East, James posed the analogy of cake cutting to the author. Among numerous relevant topics within the area of fair division, there are several types of problems that have been studied, for example: fair cake-cutting, fair chore division, fair item assignment and fair resource allocation.

The literature for fair cake-cutting dates back to at least 1948, for example see Steinhaus \cite{Steinhauscakecutting} which begins in the first paragraph by suggesting the custom ``of dividing an object into two equal parts by letting one partner halve it and the other choose his half" was already probably many centuries old. Hundreds of papers on the topic have appeared since, referencing them here would drown out the references more relevant to the bulk of this paper. 

For the point the author wishes to make we need not complicate things with who is cutting/dividing, with who is choosing, or with cakes that have different toppings. Rather than keeping the cake analogy, we equivalently consider dividing a $2$-dimensional shape in to $n \in \mathbb{Z}^+$ colours. 

\begin{definition}
	A \textit{division} of a $2$-dimensional shape in to $n \in \mathbb{Z}^+$ colours is a partition of the shape in to $n$ regions, each region coloured with a unique colour. 
\end{definition}

First let us define automorphisms for a shape division.

\begin{definition}
	An \textit{automorphism} of a shape division in to $n$ colours is any combination of rotating and/or reflecting the divided/coloured shape that leads to the exact same shape with the same division. 
\end{definition}

Note that:
\begin{enumerate}
	\item While we do not require the permuted shape to preserve colours for each region, if one region with colour $X$ is permuted to a region of colour $Y$, all regions with colour $X$ are required to permute to regions of colour $Y$; 
	\item Associated with each automorphism of a shape division is a permutation of the colours;
	\item The identity automorphism is to not rotate or reflect the shape, ie. leave it alone; and
	\item We refer to the automorphism group of a shape division as non-trivial when it does not consist of just the identity automorphism.
\end{enumerate}

\begin{definition}
	We shall refer to the division of a $2$-dimensional shape in to $n \in \mathbb{Z}^+$ colours as:
	\begin{enumerate}
		\item \textit{fair} if each colour/region fills the same area;
		\item \textit{symmetric} if the automorphism group of the shape division is non-trivial; and
		\item \textit{strongly symmetric} if the shape division has an automorphism that is not the identity but does preserve the colours of permuted regions, ie. the associated colour permutation is the identity. 
	\end{enumerate}
\end{definition}

Note that all strongly symmetric shape divisions are symmetric. Most people would agree that our definitions for fair and (strongly) symmetric in the context of shape divisions are fairly reasonable, for example:
\begin{enumerate}
	\item Our definition of symmetric covers when the colours are not really relevant, merely help to distinguish between the different regions of the division; and
	\item Our definition of strongly symmetric covers when the colours are important and we want automorphisms to preserve the colour of regions.
\end{enumerate}

It is not difficult to find examples which establish that neither fair nor (strongly) symmetric implies the other when dividing shapes in to colours, see for example Figure \ref{fig:cakecuts}.

\begin{figure}[!ht]
	\begin{center}
		\begin{tikzpicture}
			\draw[thick] (-7,4)--(-7,0)--(-1,0)--(-1,4)--(-7,4);
			\draw[thick] (-6,4)--(-6,0);
			\draw[thick] (-2,4)--(-2,0);
			\fill[pattern=north west lines]  (-7,4)--(-7,0)--(-6,0)--(-6,4);
			\fill[pattern=north west lines]  (-2,4)--(-2,0)--(-1,0)--(-1,4);
			\draw[<->,thick] (-7,4.25)--(-6,4.25);
			\node () at (-6.5,4.5) {$x$};
			\draw[<->,thick] (-6,4.25)--(-2,4.25);
			\node () at (-4,4.5) {$4x$};
			\draw[<->,thick] (-2,4.25)--(-1,4.25);
			\node () at (-1.5,4.5) {$x$};
			
			\draw[thick] (7,4)--(7,0)--(1,0)--(1,4)--(7,4);
			\draw[thick] (2.5,4)--(2.5,0);
			\draw[thick] (4,4)--(4,2);
			\draw[thick] (4,2)--(7,2);
			\fill[pattern=north west lines]  (7,4)--(7,2)--(4,2)--(4,4);
			\fill[pattern=north west lines]  (2.5,4)--(2.5,0)--(1,0)--(1,4);
			\draw[<->,thick] (7,4.25)--(4,4.25);
			\node () at (5.5,4.5) {$2x$};
			\draw[<->,thick] (4,4.25)--(2.5,4.25);
			\node () at (3.25,4.5) {$x$};
			\draw[<->,thick] (2.5,4.25)--(1,4.25);
			\node () at (1.75,4.5) {$x$};
			\draw[<->,thick] (7.25,4)--(7.25,2);
			\node () at (7.5,3) {$y$};
			\draw[<->,thick] (7.25,2)--(7.25,0	);
			\node () at (7.5,1) {$y$};
		\end{tikzpicture}
		\caption{Dividing a rectangle in to two colours/shades such that it is: (strongly) symmetric while not fair (left); and fair while not (strongly) symmetric (right).}
		\label{fig:cakecuts}
		\vspace{-0.5cm}
	\end{center}
\end{figure}	

Note that even if the rectangles from Figure \ref{fig:cakecuts} are squares, we still have that neither fair nor (strongly) symmetric implies the other. It is fairly easy to generalise the examples in Figure \ref{fig:cakecuts} with the same results for dividing a rectangle, including the case of a square, in to $n > 2$ colours. Hence, in the context of shape division, the terms fair and symmetric are far from being equivalent, and hence definitely should not be considered synonymous, to one another. A philosophical discussion could be had on whether two precisely defined terms could and/or should be considered synonymous if they are equivalent in the sense of capturing precisely the same objects.

It could be an interesting direction of research to see where fairness and (strong) symmetry do and do not overlap when dividing various different shapes into various numbers of colours. Nevertheless, it seems reasonable at this point to suggest that in a general context we should not treat the terms fair and symmetric as equivalent or synonymous to one another, as we have just seen two reasonable examples where that would break down.

Without further examination or further discussion this may leave one doubting whether notions of symmetry for games should be referred to as notions of fairness. However, astute readers may have noticed that there are some fundamental differences between our notions of symmetry for shape division compared to our notions of symmetry for games. Our symmetry notions for games have required at the very least that players be indifferent between rearrangements of positions for some transitive subgroup of the player or game permutations. Whereas:
\begin{enumerate}
	\item symmetric in the context of shape division is analogous to defining a game as symmetric if the game has an automorphism not equal to the identity; and
	\item strongly symmetric in the context of shape division is analogous to defining a game as strongly symmetric if the game has an automorphism not equal to the identity but that does use the identity player pemutation.
\end{enumerate}

Neither of these have been considered notions of symmetry for games in this paper, or really anywhere in the literature, with the exception of \cite{tohme2019structural} who define partial symmetries in games for when some players are indifferent between various positions in a game, but not indifferent between all positions. 

The notions of symmetry for games that we have considered would be closer to the following notions of symmetry for dividing shapes. 

\begin{definition}
	We shall refer to the division of a $2$-dimensional shape in to $n \in \mathbb{Z}^+$ colours as:
	\begin{enumerate}
		\item \textit{transitively symmetric} if every colour permutation in a transitive subgroup of all colour permutations is associated to at least one automorphism of the shape division; and
		\item \textit{$n$-transitively symmetric} if every colour permutation is associated to at least one automorphism of the shape division.
	\end{enumerate}
\end{definition}

Note all $n$-transitively symmetric shape divisions are transitively symmetric, and all transitively symmetric cake divisions are:
\begin{enumerate}
	\item fair;
	\item symmetric; and
	\item not necessarily strongly symmetric;
\end{enumerate}

Since both all $n$-transitively symmetric shape divisions and all transitively symmetric shape divisions are fair, it is reasonable to refer to these as notions of fairness for dividing shapes. 

Recall from Section \ref{sec:intro} that the term fair has appeared in the context of zero-sum games, noting that zero-sum games are a subclass of strategic-form games. For example a $2$-player zero-sum game is defined as fair when the (minimax or expected) value of the game is $0$, ie. when the expected payoff under perfect play for both players is $0$. This gives us another notion for zero-sum games of players being indifferent between which position they play. 

A $2$-player zero-sum game can be fair without both players having the same number of strategies, so clearly a $2$-player zero-sum game being fair does not imply that it satisfies any of the notions of symmetry defined in this paper. Nor does it really seem to make much sense to refer to a fair $2$-player zero-sum game as symmetric. 

\begin{theorem}
	If a two person zero-sum game $\Gamma$ is symmetric then it is fair (as mentioned earlier, this result is also in von Neumann and Morgenstern \cite[Pages 165-166]{VNM}).
	\begin{proof}
		Let $\pi = (12)$. It follows from $\pi \in \Aut(\Gamma)$ that:
		
		 $\displaystyle\min_{s_2}\max_{s_1}u_1(s_1, s_2) = \min_{s_2}\max_{s_1}u_{\pi(1)}(s_{\pi^{-1}(1)}, s_{\pi^{-1}(2)}) = \min_{s_2}\max_{s_1}u_2(s_2, s_1) = \min_{s_1}\max_{s_2}u_2(s_1, s_2)$.
	\end{proof}
\end{theorem}

One of our notions of fairness for games is that players be indifferent between which position they play. Game isomorphisms establish that players are indifferent between playing the mapped positions for each game, in the case of game automorphisms we get players being indifferent between playing different positions in the same game. If the automorphism group of a game is player transitive, the players are indifferent between which position they play, consequently any notion of symmetry requiring the automorphism group be player transitive falls inside the notion of fairness that players be indifferent between which position they play. 

Note however that the players may still care about the positions of their opponents, which is the stronger notion of fairness for games that we have considered.

\subsection{Classifying A Game} \label{subsec:classifying}
While our distinct symmetry notions give us various descriptive definitions of strategic fairness, they do not give us a constructive way to determine where a particular game lies. We now discuss various strategies for classifying a game which will be crucial later on when identifying examples for each combination of symmetry notions considered in this paper.

The strategies for classifying a game introduced in this subsection are not an attempt to outline all the steps required for algorithms that can be implemented to classify games, though with some gaps filled in many of the strategies could be used in such algorithms. It would be a useful future research direction to examine the complexity of deciding whether a game satisfies each notion of symmetry along with outlining algorithms for classifying a game, a potential application of doing so will be mentioned in Subsection \ref{sec:parameterisedgames}. 

To test whether a game $\Gamma$ is fully or standard symmetric: we first try to construct a matching $M$ of the strategy sets where for each profile $s \in M$, all players have the same payoff. If no such matching exists $\Gamma$ is neither fully nor standard symmetric. For example in Matching Pennies, since there does not exist a strategy profile where all players receive the same payoff, we can conclude Matching Pennies is non-standard symmetric.

If such matchings exist: to test for full symmetry we check whether such a matching induces automorphisms for permutations that generate $S_N$; and to test for standard symmetry we check whether such a matching induces automorphisms for player permutations that generate a transitive subgroup of $S_N$, noting that to conclude non-standard symmetry we must check that the game is not invariant under the bijections induced by any such matching and transitive subgroup of $S_N$. 

The reader should note that every $n$-cycle generates a transitive subgroup of $S_N$, but not all transitive subgroups of $S_N$ contain an $n$-cycle. For example the Klein group $\{e, (12) \circ (34), (13) \circ (24), (14) \circ (23)\}$ is a transitive subgroup of $S_4$ that does not contain any $4$-cycles.

To test for $n$-transitivity we check whether there exists automorphisms for permutations that generate $S_N$; and to test for symmetry (ie. transitivity) we check whether there exists automorphisms for permutations that generate a transitive subgroup of $S_N$, again noting that to conclude that a game is not symmetric we must check that the game is not invariant under any transitive subgroup of $S_{\Gamma}$. 

If we know a game is symmetric (ie. transitive) and want to show it is only-transitive, a sufficient condition is to find a strategy profile $s \in A$ whose payoffs do not appear elsewhere under all possible permutations. For example consider Example \ref{stdsymeg} and suppose it has an automorphism whose player permutation is $(23)$. The payoffs for the profile $(a, a, b)$ are $(3, 7, 4)$, so we would need a strategy profile $s \in A$ with payoffs $(3, 4, 7)$. However no such profile exists, hence Example \ref{stdsymeg} is an only-transitive standard symmetric game.

\subsection{Parameterised Symmetric Games} \label{sec:parameterisedgames}
Given a subset $G$ of game bijections we construct the \textit{parameterised game $\Gamma(G)$ of $G$} as follows: for each $g \in \langle{G}\rangle$, $s \in A$ and $i \in N$, set $u_i(s) = u_{g(i)}\bigl(g(s)\bigr)$. Since automorphisms are closed under composition we have $\langle{G}\rangle \leq \Aut(\Gamma)$, hence each orbit of $(N\times{A})/\langle{G}\rangle$ has the same payoff. An algorithm for this construction method has been implemented by the author in C++, the code is available at \cite{GLCode}.

\begin{example} \label{constructeg}
	Let $G = \{\bigl((12) ; \bigl(\begin{smallmatrix} a & b \\ c & d \end{smallmatrix}\bigr), \bigl(\begin{smallmatrix} c & d \\ a & b \end{smallmatrix}\bigr)\bigr)\}$. For $\Gamma(G)$ we require:
	\begin{align*}
		u_1(a, c) &= u_2(a, c) = \alpha   &   u_1(a, d) &= u_2(b, c) = \gamma \\
		u_1(b, c) &= u_2(a, d) = \beta   &   u_1(b, d) &= u_2(b, d) = \delta
	\end{align*} 
	\begin{center}
	\begin{game}{2}{2}[$\Gamma(G)$]
			  \> $c$  \> $d$ \\
		$a$   \> $\alpha, \alpha $  \> $\gamma, \beta $ \\
		$b$   \> $\beta, \gamma $  \> $\delta, \delta $
	\end{game}
	\end{center}
	
	We call $\alpha, \beta, \gamma, \delta \in \mathbb{R}$ the \textit{parameters} of $\Gamma(G)$. Note that distinct parameter choices may lead to strategically inequivalent games, even though both games will have the same automorphism group. All fully symmetric $2$-player $2$-strategy games are isomorphic to $\Gamma(G)$ for at least one choice of parameters, hence $\Gamma(G)$ is a general form for fully symmetric $2$-player $2$-strategy games, or equivalently standard symmetric $2$-player $2$-strategy games.
\end{example}

	We can define a partial order $\leq$ on parameterised games as follows: $\Gamma(G) \leq \Gamma(G')$ when given a set of parameter choices for $\Gamma(G')$ there exists a set of parameter choices for $\Gamma(G)$ such that $\Gamma(G) \cong \Gamma(G')$. We illustrate our order in Figures \ref{2pHasse} and \ref{3pHasse} using the Hasse diagrams for $\leq$ on parameterised symmetric $2$-player and $3$-player $2$-strategy games up to isomorphism, which were constructed using the code at \cite{GLCode}.

	\begin{figure}[!ht]
		\begin{center}
		\begin{tikzpicture}
    			\node (r3c1) at (0,0)  {\begin{tabular}{| c | c |}
    			\multicolumn{2}{c}{$\Gamma(G_{31})$}\\
    			\hline
    			$\alpha, \alpha$  & $\alpha, \alpha$ \\ \hline
			$\alpha, \alpha$  & $\alpha, \alpha$ \\
    			\hline
  		\end{tabular}};
	
    			\node (r2c2) at (1.5,-2) {\begin{tabular}{| c | c |}
    			\multicolumn{2}{c}{$\Gamma(G_{22})$}\\
    			\hline
    			$\alpha, \beta$  & $\beta, \alpha$ \\ \hline
			$\beta, \alpha$  & $\alpha, \beta$ \\
    			\hline
  		\end{tabular}};	
	
			\node (r2c1) at (-1.5,-2)  {\begin{tabular}{| c | c |}
    			\multicolumn{2}{c}{$\Gamma(G_{21})$}\\
    			\hline
    			$\alpha, \alpha$  & $\beta, \beta$ \\ \hline
			$\beta, \beta$  & $\alpha, \alpha$ \\
    			\hline
  		\end{tabular}};

			\node (r1c1) at (-1.5,-4)  {\begin{tabular}{| c | c |}
    			\multicolumn{2}{c}{$\Gamma(G_{11})$}\\
    			\hline
    			$\alpha, \alpha$  & $\beta, \gamma$ \\ \hline
			$\gamma, \beta$  & $\delta, \delta$ \\
    			\hline
  		\end{tabular}};
  		
  		\draw (r2c1) -- (r3c1) -- (r2c2);
  		\draw (r2c1) -- (r1c1);
		\end{tikzpicture}
		\end{center}
		
		\begin{center}
		$G_{11} = \{\bigl((12) ; \bigl(\begin{smallmatrix} a & b \\ c & d \end{smallmatrix}\bigr), \bigl(\begin{smallmatrix} c & d \\ a & b \end{smallmatrix}\bigr)\bigr)\}$, $G_{21} = G_{11} \cup \{\bigl((12) ; \bigl(\begin{smallmatrix} a & b \\ d & c \end{smallmatrix}\bigr), \bigl(\begin{smallmatrix} c & d \\ b & a \end{smallmatrix}\bigr)\bigr)\}$,
		
		$G_{22} = \{\bigl((12) ; \bigl(\begin{smallmatrix} a & b \\ d & c \end{smallmatrix}\bigr), \bigl(\begin{smallmatrix} c & d \\ a & b \end{smallmatrix}\bigr)\bigr)\}$, $G_{31} = G_{11} \cup G_{22}$.
		\end{center}
		\vspace{-0.5cm}
		\caption{A plot of the graph of the Hasse diagram for $\leq$ on parameterised symmetric $2$-player $2$-strategy games up to isomorphism.}
		\label{2pHasse}
	\end{figure}
	
	One future direction of research would be to compute the Hasse diagrams for all (symmetric) parameterised games up to isomorphism for a fixed number of players and fixed number of strategies, seeing which combinations of player and strategy counts are computationally feasible with modern day hardware and compilers. One can then also get the poset of (symmetric) games up to isomorphism by finding the strategically inequivalent games for each parameterised game. It would be good to come up with a way to allow us to confirm/verify the result from Goforth and Robinson \cite{GoforthRobinson} that there are 144 ordinal equivalence classes for the 2-player 2-strategy games, which should be possible using game bijections, and obtain numbers for the number of (parameterised) (symmetric) games up to isomorphism for various player and strategy counts.
	
	In order to achieve the above research goals, a number of algorithms beyond those at \cite{GLCode} would need to be examined and implemented, including:
	\begin{enumerate}
		\item Implement algorithms to classify a game for each desired symmetry notion, see Subsection \ref{subsec:classifying} for a discussion on this;
		\item Find a precise and accurate definition for isomorphisms between two parameterised games, which is likely more complicated than for non-parameterised games due to the utilities for outcomes/profiles no longer being ordered;
		\item Implement an algorithm to check whether there exists an isomorphism between two parameterised games, ie. check whether two parameterised games are equivalent. This may require iterating through every relabelling of the players and strategies;
		\item Implement an algorithm for the partial-order we defined on parameterised games;
		\item Implement a struct for posets of strategic-form games ordered using the partial-order algorithm; and
		\item Unless the search space can be reduced, iterate through all appropriate subsets/subgroups of the game bijections, constructing the game for each subset/subgroup, checking that it meets the desired properties in the case of symmetric games, checking whether each new constructed game is equivalent to any games already found, and if it not equivalent to any games already found then insert it in to the poset. 
	\end{enumerate}
	
	When the above is successfully achieved, one will also need a way to output the Hasse diagrams in a way that humans can interpret, preferably also outputting the Hasse diagrams as Ti\emph{K}Z code or code for any other {\LaTeX} package.
	
	The author sees little point in constructing the Hasse diagram for non-parameterised games as it will essentially be the same as the Hasse diagram for parameterised games except each parameterised game is replaced with an equivalence class of strategically inequivalent games that essentially have the same automorphism group.
	
	Subsection \ref{subsec:examples} contains examples of games for the remaining combinations of symmetry notions that examples have not already been given for. These examples were constructed using the code at \cite{GLCode} and using the following strategies.
    
	To construct a symmetric game or an $n$-transitively symmetric game we use bijections that generate a player transitive  or player $n$-transitive subgroup respectively. 
	
	To construct an only-transitive symmetric game it is not sufficient to use bijections that generate an only-transitive subgroup, we must construct $\Gamma(G)$ and check that it is only-transitive. This is due to $\langle{G}\rangle$ possibly being a proper subgroup of $\Aut(\Gamma)$. For example, if we take:
	\[G = \{ \bigl((123) ; \bigl(\begin{smallmatrix} a & b \\ d & c \end{smallmatrix}\bigr), \bigl(\begin{smallmatrix} c & d \\ e & f \end{smallmatrix}\bigr), \bigl(\begin{smallmatrix} e & f \\ a & b \end{smallmatrix}\bigr)\bigr),
               \bigl((123) ; \bigl(\begin{smallmatrix} a & b \\ c & d \end{smallmatrix}\bigr), \bigl(\begin{smallmatrix} c & d \\ f & e \end{smallmatrix}\bigr), \bigl(\begin{smallmatrix} e & f \\ a & b \end{smallmatrix}\bigr)\bigr) \},\] 
    then $N\times{A}$ has one orbit under $\langle{G}\rangle$ (i.e. $\Gamma(G)$ has one parameter/payoff) despite $\langle{G}\rangle$ being an only-transitive subgroup.    
        
	To construct a standard symmetric game we use the bijections induced from a matching of the strategy sets and player permutations which generate a transitive subgroup of $S_N$. To construct a non-standard symmetric game, we first choose game bijections which are not obviously from the same matching, construct $\Gamma(G)$ and check whether it is non-standard symmetric. We construct fully and non-fully symmetric games similarly.
	
\subsection{Further examples} \label{subsec:examples}
So far we have seen examples of fully symmetric, only-transitive standard symmetric and $n$-transitively non-standard symmetric games. We now look at examples constructed with the code at \cite{GLCode} to show that our notions of symmetry are related as shown in the Euler diagram in Figure \ref{fig:Eulerdiag}. A similar approach with regards to constructing examples has been taken by the author and East in \cite{EastHam} within the context of lattice path enumeration to give the precise relationships between finiteness properties of path counts to end points, geometrical properties of the underlying step set, algebraic properties of the monoid of end points, and combinatorial properties of a certain bi-labelled digraph naturally associated to the underlying step set. Knowing the precise relationship between the various notions of symmetry is useful for:
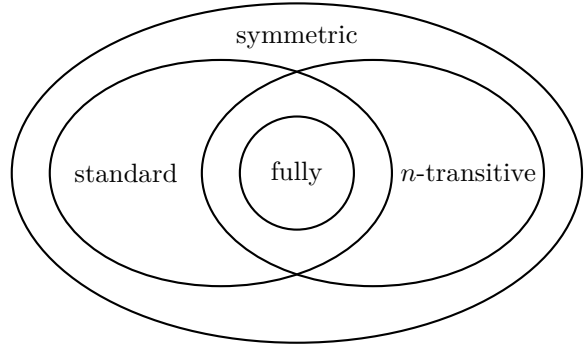
\begin{enumerate}
	\item The theory of (symmetric) games and more generally from a pure mathematics perspective; and
	\item Identifying which combinations of symmetry notions may appear in areas like artificial intelligence, biology, computer science, economics, legal systems, logic, philosophy, politics, along with social choice and voting theory. Though a more thorough investigation of examples for each feasible combination would be needed to identify which strategic situations arise in different contexts. 
\end{enumerate}

\begin{figure}[!ht]
	\begin{center}
	\begin{tikzpicture}
		nodes={draw, ultra thick}
		\draw[thick] (0,0) circle (0.75cm);			\draw (0, 0) node {fully};
		\draw[thick] (-1.5,0) ellipse (2.75cm and 1.5cm);	\draw (-2.5,0) node {standard};
		\draw[thick] (1.5,0) ellipse (2.75cm and 1.5cm);	\draw (2.5,0) node {$n$-transitive};
		\draw[thick] (0,0) ellipse (5.5cm and 2.25cm);	\draw (0,1.8) node {symmetric};
	\end{tikzpicture}
	\end{center}
	\vspace{-0.5cm}
	\caption{Euler diagram for label-independent symmetry notions.}
	\label{fig:Eulerdiag}
\end{figure}

\begin{example} \label{egntransstandnonfull}	
	An $n$-transitively non-fully standard symmetric $3$-player game.
	\begin{center}
	\begin{game}{2}{2}[$(a,,)$]
			\>  $e$      \>  $f$      \\
		$c$ \>  $\alpha,\alpha,\alpha$  \>  $\beta,\gamma,\delta$  \\
		$d$ \>  $\gamma,\delta,\beta$  \>  $\delta,\gamma,\beta$  
	\end{game}
	\hspace*{5mm}
	\begin{game}{2}{2}[$(b,,)$]
			\>  $e$      \>  $f$      \\
		$c$ \>  $\delta,\beta,\gamma$  \>  $\beta,\delta,\gamma$  \\
		$d$ \>  $\gamma,\beta,\delta$  \>  $\alpha,\alpha,\alpha$  
	\end{game}
	\end{center}
	
	\begin{align*}
		G = \{\bigl((123) ; \bigl(\begin{smallmatrix} a & b \\ c & d \end{smallmatrix}\bigr), \bigl(\begin{smallmatrix} c & d \\ e & f \end{smallmatrix}\bigr), \bigl(\begin{smallmatrix} e & f \\ a & b \end{smallmatrix}\bigr)\bigr), \bigl((12) ; \bigl(\begin{smallmatrix} a & b \\ d & c \end{smallmatrix}\bigr), \bigl(\begin{smallmatrix} c & d \\ b & a \end{smallmatrix}\bigr), \bigl(\begin{smallmatrix} e & f \\ f & e \end{smallmatrix}\bigr)\bigr)\}
	\end{align*}	

    Since $\langle{G}\rangle$ is $n$-transitive, and the first generator generates a player transitive and strategy trivial group with the matching $M  = \{(a,c,e), (b,d,f)\}$, $\Gamma(G)$ is $n$-transitively and standard symmetric. Furthermore since the bijections induced by $M$ from player transpositions are not automorphisms, $\Gamma(G)$ is non-fully symmetric.
\end{example}

\begin{figure} 
		\begin{center}
		\begin{tikzpicture}
    			\node (r4c1) at (0,0)  {\begin{tabular}{| c | c |}
    			\multicolumn{2}{c}{$\Gamma(G_{41})$}\\
    			\hline
    			$\alpha,\alpha,\alpha$ & $\alpha,\alpha,\alpha$ \\ \hline
				$\alpha,\alpha,\alpha$ & $\alpha,\alpha,\alpha$ \\
    			\hline \hline
				$\alpha,\alpha,\alpha$ & $\alpha,\alpha,\alpha$ \\ \hline
				$\alpha,\alpha,\alpha$ & $\alpha,\alpha,\alpha$ \\
    			\hline
  		\end{tabular}};
			
    			\node (r3c1)  at (-4,-3)    {\begin{tabular}{| c | c |}
    			\multicolumn{2}{c}{$\Gamma(G_{31})$}\\
    			\hline
    			$\alpha,\alpha,\alpha$ & $\beta,\beta,\beta$ \\ \hline
			$\beta,\beta,\beta$ & $\alpha,\alpha,\alpha$ \\
    			\hline \hline
			$\beta,\beta,\beta$ & $\alpha,\alpha,\alpha$ \\ \hline
			$\alpha,\alpha,\alpha$ & $\beta,\beta,\beta$ \\
    			\hline
  		\end{tabular}};
			
			\node (r3c2)  at (0,-3)     {\begin{tabular}{| c | c |}
			\multicolumn{2}{c}{$\Gamma(G_{32})$}\\
    			\hline
    			$\alpha,\alpha,\alpha$ & $\beta,\beta,\delta$ \\ \hline
			$\beta,\delta,\beta$ & $\delta,\beta,\beta$ \\
    			\hline \hline
			$\delta,\beta,\beta$ & $\beta,\delta,\beta$ \\ \hline
			$\beta,\beta,\delta$ & $\alpha,\alpha,\alpha$  \\
    			\hline
  		\end{tabular}};
			
			\node (r2c1)  at (-4,-6)    {\begin{tabular}{| c | c |}
			\multicolumn{2}{c}{$\Gamma(G_{21})$}\\
    			\hline
    			$\alpha,\alpha,\alpha$ & $\beta,\beta,\delta$ \\ \hline
			$\beta,\delta,\beta$ & $\sigma,\rho,\rho$ \\
    			\hline \hline
			$\delta,\beta,\beta$ & $\rho,\sigma,\rho$ \\ \hline
			$\rho,\rho,\sigma$ & $\omega,\omega,\omega$  \\
    			\hline
  		\end{tabular}};
			
			\node (r2c2)  at  (0,-6)    {\begin{tabular}{| c | c |}
			\multicolumn{2}{c}{$\Gamma(G_{22})$}\\
    			\hline
    			$\alpha,\alpha,\alpha$ & $\beta,\gamma,\delta$ \\ \hline
			$\gamma,\delta,\beta$ & $\delta,\gamma,\beta$ \\
    			\hline \hline
			$\delta,\beta,\gamma$ & $\beta,\delta,\gamma$ \\ \hline
			$\gamma,\beta,\delta$ & $\alpha,\alpha,\alpha$  \\
    			\hline
  		\end{tabular}};
			
			\node (r2c3)  at (4,-6)     {\begin{tabular}{| c | c |}
			\multicolumn{2}{c}{$\Gamma(G_{23})$}\\
    			\hline
    			$\alpha,\alpha,\alpha$ & $\beta,\gamma,\delta$ \\ \hline
			$\gamma,\delta,\beta$ & $\delta,\beta,\gamma$ \\
    			\hline \hline
			$\delta,\beta,\gamma$ & $\gamma,\delta,\beta$ \\ \hline
			$\beta,\gamma,\delta$ & $\alpha,\alpha,\alpha$  \\
    			\hline
  		\end{tabular}};	
			
			\node (r1c1)  at (0,-9)     {\begin{tabular}{| c | c |}
			\multicolumn{2}{c}{$\Gamma(G_{11})$}\\
    			\hline
    			$\alpha,\alpha,\alpha$ & $\beta,\gamma,\delta$ \\ \hline
			$\gamma,\delta,\beta$ & $\sigma,\rho,\tau$ \\
    			\hline \hline
			$\delta,\beta,\gamma$ & $\tau,\sigma,\rho$ \\ \hline
			$\rho,\tau,\sigma$ & $\omega,\omega,\omega$  \\
    			\hline
  		\end{tabular}};
  			
  		\draw (r3c1) -- (r4c1) -- (r3c2);
  		\draw (r3c1) -- (r2c1) -- (r3c2) -- (r2c2);
  		\draw (r3c2) -- (r2c3);
  		\draw (r2c1) -- (r1c1) -- (r2c2);
  		\draw (r2c3) -- (r1c1);
			
		\end{tikzpicture}
		\end{center}
		
		\begin{center}
		$G_{11} = \{\bigl((123) ; \bigl(\begin{smallmatrix} a & b \\ c & d \end{smallmatrix}\bigr), \bigl(\begin{smallmatrix} c & d \\ e & f \end{smallmatrix}\bigr), \bigl(\begin{smallmatrix} e & f \\ a & b \end{smallmatrix}\bigr)\bigr)\}$, 
		
		$G_{21} = G_{11} \cup \{\bigl((12) ; \bigl(\begin{smallmatrix} a & b \\ c & d \end{smallmatrix}\bigr), \bigl(\begin{smallmatrix} c & d \\ a & b \end{smallmatrix}\bigr), \bigl(\begin{smallmatrix} e & f \\ e & f \end{smallmatrix}\bigr)\bigr)\}$,
		
		$G_{22} = G_{11} \cup \{\bigl((12) ; \bigl(\begin{smallmatrix} a & b \\ d & c \end{smallmatrix}\bigr), \bigl(\begin{smallmatrix} c & d \\ b & a \end{smallmatrix}\bigr), \bigl(\begin{smallmatrix} e & f \\ f & e \end{smallmatrix}\bigr)\bigr)\}$,
		
		$G_{23} = \{\bigl((123) ; \bigl(\begin{smallmatrix} a & b \\ d & c \end{smallmatrix}\bigr), \bigl(\begin{smallmatrix} c & d \\ f & e \end{smallmatrix}\bigr), \bigl(\begin{smallmatrix} e & f \\ b & a \end{smallmatrix}\bigr)\bigr)\}$,
		
		$G_{31} = G_{21} \cup \{\bigl((123) ; \bigl(\begin{smallmatrix} a & b \\ d & c \end{smallmatrix}\bigr), \bigl(\begin{smallmatrix} c & d \\ f & e \end{smallmatrix}\bigr), \bigl(\begin{smallmatrix} e & f \\ a & b \end{smallmatrix}\bigr)\bigr)\}$,
		
		$G_{32} = G_{2i} \cup G_{2j}$ for all distinct $i, j \in \{1,2,3\}$,
		
		$G_{41} = G_{31} \cup G_{32}$.
		\end{center}
		
		\vspace{-0.5cm}
		\caption{A plot of the graph of the Hasse diagram for $\leq$ on parameterised symmetric $3$-player $2$-strategy games up to isomorphism.}
		\label{3pHasse}
	\end{figure}
	
Cheng et al. \cite{CRVWSym} showed that fully symmetric $2$-strategy games have at least one pure strategy Nash equilibrium. They also noted that Rock, Paper, Scissors is an example of a fully symmetric $2$-player $3$-strategy game with no pure strategy Nash equilibria, and indirectly that Matching Pennies is an example of a non-standard symmetric $2$-player $2$-strategy game which has no pure strategy Nash equilibria. The reader may like to verify that Example \ref{stdsymeg} is a standard symmetric $2$-strategy game with no pure strategy Nash equilibria. 

Note Example \ref{egntransstandnonfull} is the only parameterised $n$-transitively non-fully standard symmetric $3$-player $2$-strategy game up to isomorphism. Furthermore note there are pure strategy Nash equilibria for each choice of parameters. The author has been unable to show whether the result from Cheng et al. \cite{CRVWSym} weakens to $n$-transitively standard symmetric $2$-strategy games.

\begin{example} \label{OTNSeg}
Two only-transitive non-standard symmetric $4$-player games.
\begin{center}
  \begin{game}{2}{2}[$(a,c,,)$]
       \>  $g$          \>  $h$          \\
    $e$  \>  $\alpha,\beta,\gamma,\delta$  \>  $\rho,\tau,\sigma,\omega$  \\
    $f$  \>  $\sigma,\omega,\rho,\tau$  \>  $\omega,\rho,\tau,\sigma$  
  \end{game}
  \hspace*{5mm}
  \begin{game}{2}{2}[$(a,d,,)$]
       \>  $g$          \>  $h$          \\
    $e$  \>  $\delta,\alpha,\beta,\gamma$  \>  $\tau,\sigma,\omega,\rho$  \\
    $f$  \>  $\gamma,\delta,\alpha,\beta$  \>  $\beta,\gamma,\delta,\alpha$  
  \end{game}
  \\
  \begin{game}{2}{2}[$(b,c,,)$]
       \>  $g$          \>  $h$          \\
    $e$  \>  $\beta,\gamma,\delta,\alpha$  \>  $\gamma,\delta,\alpha,\beta$  \\
    $f$  \>  $\tau,\sigma,\omega,\rho$  \>  $\delta,\alpha,\beta,\gamma$  
  \end{game}
  \hspace*{5mm}
  \begin{game}{2}{2}[$(b,d,,)$]
       \>  $g$          \>  $h$          \\
    $e$  \>  $\omega,\rho,\tau,\sigma$  \>  $\sigma,\omega,\rho,\tau$  \\
    $f$  \>  $\rho,\tau,\sigma,\omega$  \>  $\alpha,\beta,\gamma,\delta$  
  \end{game}
\end{center}

\[ G = \{ \bigl((1234) ; \bigl(\begin{smallmatrix} a & b \\ d & c \end{smallmatrix}\bigr), \bigl(\begin{smallmatrix} c & d \\ e & f \end{smallmatrix}\bigr), \bigl(\begin{smallmatrix} e & f \\ g & h \end{smallmatrix}\bigr), \bigl(\begin{smallmatrix} g & h \\ a & b \end{smallmatrix}\bigr)\bigr)\} \]

Since there does not exist any profile where the payoffs are equal and $\langle{G}\rangle$ is transitive, $\Gamma(G)$ is non-standard symmetric. Now for the strategy profile $(a,c,e,g)$ we have payoffs $(\alpha, \beta, \gamma, \delta)$. If $\Gamma(G)$ had an automorphism using $(23)$ then there would be a strategy profile $s \in A$ with payoffs $(\alpha, \gamma, \beta, \delta)$. Since no such profile exists $\Gamma(G)$ is only-transitive.

\begin{center}
  \begin{game}{2}{2}[$(a,c,,)$]
       \>  $g$          \>  $h$          \\
    $e$  \>  $\alpha,\alpha,\beta,\beta$  \>  $\gamma,\delta,\delta,\gamma$  \\
    $f$  \>  $\delta,\gamma,\gamma,\delta$  \>  $\beta,\beta,\alpha,\alpha$  
  \end{game}
  \hspace*{5mm}
  \begin{game}{2}{2}[$(a,d,,)$]
       \>  $g$          \>  $h$          \\
    $e$  \>  $\gamma,\delta,\delta,\gamma$  \>  $\alpha,\alpha,\beta,\beta$  \\
    $f$  \>  $\beta,\beta,\alpha,\alpha$  \>  $\delta,\gamma,\gamma,\delta$  
  \end{game}
  \\
  \begin{game}{2}{2}[$(b,c,,)$]
       \>  $g$          \>  $h$          \\
    $e$  \>  $\delta,\gamma,\gamma,\delta$  \>  $\beta,\beta,\alpha,\alpha$  \\
    $f$  \>  $\alpha,\alpha,\beta,\beta$  \>  $\gamma,\delta,\delta,\gamma$  
  \end{game}
  \hspace*{5mm}
  \begin{game}{2}{2}[$(b,d,,)$]
       \>  $g$          \>  $h$          \\
    $e$  \>  $\beta,\beta,\alpha,\alpha$  \>  $\delta,\gamma,\gamma,\delta$  \\
    $f$  \>  $\gamma,\delta,\delta,\gamma$  \>  $\alpha,\alpha,\beta,\beta$  
  \end{game}
\end{center}

\begin{align*}
G' = \{ \bigl((12) \circ (34) &; \bigl(\begin{smallmatrix} a & b \\ d & c \end{smallmatrix}\bigr), \bigl(\begin{smallmatrix} c & d \\ a & b \end{smallmatrix}\bigr), \bigl(\begin{smallmatrix} e & f \\ h & g \end{smallmatrix}\bigr), \bigl(\begin{smallmatrix} g & h \\ e & f \end{smallmatrix}\bigr)\bigr), \\
               \bigl((13) \circ (24) &; \bigl(\begin{smallmatrix} a & b \\ f & e \end{smallmatrix}\bigr), \bigl(\begin{smallmatrix} c & d \\ h & g \end{smallmatrix}\bigr), \bigl(\begin{smallmatrix} e & f \\ a & b \end{smallmatrix}\bigr), \bigl(\begin{smallmatrix} g & h \\ c & d \end{smallmatrix}\bigr)\bigr), \\
               \bigl((14) \circ (23) &; \bigl(\begin{smallmatrix} a & b \\ h & g \end{smallmatrix}\bigr), \bigl(\begin{smallmatrix} c & d \\ f & e \end{smallmatrix}\bigr), \bigl(\begin{smallmatrix} e & f \\ c & d \end{smallmatrix}\bigr), \bigl(\begin{smallmatrix} g & h \\ a & b \end{smallmatrix}\bigr)\bigr)\}
\end{align*}

That $\Gamma(G')$ is only-transitive non-standard symmetric follows by the same argument used for $\Gamma(G)$. 
\end{example}

\begin{example}
An $n$-transitively non-standard symmetric $4$-player game.
\begin{center}
\begin{game}{2}{2}[$(a,c,,)$]
       \>  $g$          \>  $h$ \\
	$e$  \>  $\alpha,\beta,\beta,\beta$  \>  $\beta,\alpha,\beta,\beta$  \\
    $f$  \>  $\beta,\beta,\beta,\alpha$  \>  $\beta,\beta,\beta,\alpha$  
\end{game}
\hspace*{5mm}
\begin{game}{2}{2}[$(a,d,,)$]
       \>  $g$          \>  $h$          \\
    $e$  \>  $\beta,\beta,\alpha,\beta$  \>  $\beta,\alpha,\beta,\beta$  \\
    $f$  \>  $\beta,\beta,\alpha,\beta$  \>  $\alpha,\beta,\beta,\beta$  
\end{game}
\\
\begin{game}{2}{2}[$(b,c,,)$]
       \>  $g$          \>  $h$          \\
    $e$  \>  $\alpha,\beta,\beta,\beta$  \>  $\beta,\beta,\alpha,\beta$  \\
    $f$  \>  $\beta,\alpha,\beta,\beta$  \>  $\beta,\beta,\alpha,\beta$  
\end{game}
\hspace*{5mm}
\begin{game}{2}{2}[$(b,d,,)$]
       \>  $g$          \>  $h$          \\
    $e$  \>  $\beta,\beta,\beta,\alpha$  \>  $\beta,\beta,\beta,\alpha$  \\
    $f$  \>  $\beta,\alpha,\beta,\beta$  \>  $\alpha,\beta,\beta,\beta$  
\end{game}
\end{center}

\[G = \{ \bigl((1234) ; \bigl(\begin{smallmatrix} a & b \\ c & d \end{smallmatrix}\bigr), \bigl(\begin{smallmatrix} c & d \\ e & f \end{smallmatrix}\bigr), \bigl(\begin{smallmatrix} e & f \\ h & g \end{smallmatrix}\bigr), \bigl(\begin{smallmatrix} g & h \\ a & b \end{smallmatrix}\bigr)\bigr), 
               \bigl((12) ; \bigl(\begin{smallmatrix} a & b \\ c & d \end{smallmatrix}\bigr), \bigl(\begin{smallmatrix} c & d \\ a & b \end{smallmatrix}\bigr), \bigl(\begin{smallmatrix} e & f \\ e & f \end{smallmatrix}\bigr), \bigl(\begin{smallmatrix} g & h \\ h & g \end{smallmatrix}\bigr)\bigr)\}\]

$\Gamma(G)$ is $n$-transitive since $\langle{G}\rangle$ is $n$-transitive, and non-standard symmetric since there does not exist any profile where all players receive the same payoff.
\end{example}

\begin{example} \label{sixplayereg} 
An only-transitive non-standard symmetric $6$-player game that has a subgroup $\langle{G}\rangle$ isomorphic to $\overrightarrow{\langle{G}\rangle}$ with $\langle{G}\rangle_N = \{\id_{\Gamma}\}$.
  
\begin{center}
\begin{game}{2}{2}[$(a,c,e,g,,)$]
	\>$k$                 \>$l$  \\
	$i$\>$1,2,1,2,1,2$     \>$3,4,5,6,7,8$\\
	$j$\>$9,10,11,12,13,14$\>$15,16,17,18,19,20$  
\end{game}
\hspace*{5mm}
\begin{game}{2}{2}[$(a,c,e,h,,)$]
	\>$k$                  \>$l$ \\
	$i$\>$5,6,7,8,3,4$      \>$20,15,19,17,18,16$\\
	$j$\>$21,22,23,24,25,26$\>$27,27,28,28,28,27$  
\end{game}
\\
\begin{game}{2}{2}[$(a,c,f,g,,)$]
	\>$k$                  \>$l$ \\
	$i$\>$11,12,13,14,9,10$ \>$29,29,30,30,30,29$\\
	$j$\>$26,24,22,23,21,25$\>$4,8,6,7,5,3$        
\end{game}
\hspace*{5mm}
\begin{game}{2}{2}[$(a,c,f,h,,)$]
	\>$k$                  \>$l$ \\
	$i$\>$17,18,19,20,15,16$\>$8,3,7,5,6,4$\\
	$j$\>$31,32,32,32,31,31$\>$16,20,18,19,17,15$  
\end{game}
\\
\begin{game}{2}{2}[$(a,d,e,g,,)$]
	\>$k$                  \>$l$   \\
	$i$\>$7,8,3,4,5,6$      \>$18,16,20,15,19,17$\\
	j\>$30,29,29,29,30,30$\>$6,4,8,3,7,5$        
\end{game}
\hspace*{5mm}
\begin{game}{2}{2}[$(a,d,e,h,,)$]
	\>$k$                  \>$l$  \\
	$i$\>$19,17,18,16,20,15$\>$32,31,32,31,32,31$\\
	$j$\>$13,11,12,10,14,9$ \>$22,26,24,25,23,21$  
\end{game}
\\
\begin{game}{2}{2}[$(a,d,f,g,,)$]
	\>$k$                  \>$l$   \\
	$i$\>$23,24,25,26,21,22$\>$12,10,14,9,13,11$\\
	$j$\>$14,12,10,11,9,13$ \>$2,2,2,1,1,1$       
\end{game}
\hspace*{5mm}
\begin{game}{2}{2}[$(a,d,f,h,,)$]
	\>$k$                  \>$l$  \\
	$i$\>$28,28,28,27,27,27$\>$24,25,23,21,22,26$\\
	$j$\>$25,23,24,22,26,21$\>$10,14,12,13,11,9$   
\end{game}
\\
\begin{game}{2}{2}[$(b,c,e,g,,)$]
	\>$k$                  \>$l$  \\
	$i$\>$13,14,9,10,11,12$ \>$25,26,21,22,23,24$\\
	$j$\>$21,25,26,24,22,23$\>$31,31,31,32,32,32$  
\end{game}
\hspace*{5mm}
\begin{game}{2}{2}[$(b,c,e,h,,)$]
	\>$k$                  \>$l$  \\
	$i$\>$30,30,30,29,29,29$\>$14,9,13,11,12,10$\\
	$j$\>$9,13,14,12,10,11$ \>$26,21,25,23,24,22$  
\end{game}
\\
\begin{game}{2}{2}[$(b,c,f,g,,)$]
	\>$k$                  \>$l$  \\
	$i$\>$22,23,21,25,26,24$\>$10,11,9,13,14,12$\\
	$j$\>$27,28,27,28,27,28$\>$16,17,15,19,20,18$  
\end{game}
\hspace*{5mm}
\begin{game}{2}{2}[$(b,c,f,h,,)$]
	\>$k$                  \>$l$  \\
	$i$\>$6,7,5,3,4,8$      \>$2,1,1,1,2,2$\\
	$j$\>$15,19,20,18,16,17$\>$4,5,3,7,8,6$  
\end{game}
\\
\begin{game}{2}{2}[$(b,d,e,g,,)$]
	\>$k$                  \>$l$  \\
	$i$\>$19,20,15,16,17,18$\>$28,27,27,27,28,28$\\
	$j$\>$5,3,4,8,6,7$      \>$17,15,16,20,18,19$  
\end{game}
\hspace*{5mm}
\begin{game}{2}{2}[$(b,d,e,h,,)$]
	\>$k$            \>$l$ \\
	$i$\>$7,5,6,4,8,3$\>$23,21,22,26,24,25$\\
	$j$\>$1,1,2,2,2,1$\>$11,9,10,14,12,13$   
\end{game}
\\
\begin{game}{2}{2}[$(b,d,f,g,,)$]
	\>$k$                  \>$l$ \\
	$i$\>$32,32,31,31,31,32$\>$24,22,26,21,25,23$\\
	$j$\>$20,18,16,17,15,19$\>$8,6,4,5,3,7$        
\end{game}
\hspace*{5mm}
\begin{game}{2}{2}[$(b,d,f,h,,)$]
	\>$k$                  \>$l$ \\
	$i$\>$18,19,17,15,16,20$\>$12,13,11,9,10,14$\\
	$j$\>$3,7,8,6,4,5$      \>$29,30,29,30,29,30$  
\end{game}
\\
\vspace{-6mm}
\end{center}

\begin{align*}
G = \{ \bigl((14) \circ (25); &\bigl(\begin{smallmatrix} a & b \\ h & g \end{smallmatrix}\bigr), \bigl(\begin{smallmatrix} c & d \\ i & j \end{smallmatrix}\bigr), \bigl(\begin{smallmatrix} e & f \\ f & e \end{smallmatrix}\bigr), \bigl(\begin{smallmatrix} g & h \\ b & a \end{smallmatrix}\bigr), \bigl(\begin{smallmatrix} i & j \\ c & d \end{smallmatrix}\bigr), \bigl(\begin{smallmatrix} k & l \\ l & k \end{smallmatrix}\bigr)\bigr),\\
               \bigl((135) \circ (246); &\bigl(\begin{smallmatrix} a & b \\ e & f \end{smallmatrix}\bigr), \bigl(\begin{smallmatrix} c & d \\ g & h \end{smallmatrix}\bigr), \bigl(\begin{smallmatrix} e & f \\ i & j \end{smallmatrix}\bigr), \bigl(\begin{smallmatrix} g & h \\ k & l \end{smallmatrix}\bigr), \bigl(\begin{smallmatrix} i & j \\ a & b \end{smallmatrix}\bigr), \bigl(\begin{smallmatrix} k & l \\ c & d \end{smallmatrix}\bigr)\bigr)\}
\end{align*}

Since there does not exist any profile where the payoffs are equal and $\langle{G}\rangle$ is transitive, $\Gamma$ is non-standard symmetric. Now the payoffs for $(a,c,e,g,i,k)$ are $(1, 2, 1, 2, 1, 2)$. If there was an automorphism for $(12)$ then there would be $s \in A$ with payoffs $(2, 1, 1, 2, 1, 2)$. Since no such profile exists $\Gamma$ is only-transitive. It can be verified that $\langle{G}\rangle$ has order $12$, which is equal to the order of $\langle (14) \circ (25), (135) \circ (246)\rangle$, hence $\langle{G}\rangle \cong \langle (14) \circ (25), (135) \circ (246)\rangle$ and $\langle{G}\rangle_N = \{\id_{\Gamma}\}$.
\end{example}

    \section*{Acknowledgements}
    		Part of the time spent conducting research for and typesetting this paper was supported by a Tasmania Graduate Research Scholarship (186). The author would especially like to express his gratitude towards:
    		\begin{enumerate}
    			\item Des FitzGerald, without whom the author would unlikely have been able to explore symmetric games for the thesis component of his honours year at the University of Tasmania, and who has been useful for advice since;
	    		\item Noah Stein who introduced the author to his doctoral work at MIT after commenting about his work on what was a discussion based website called Reddit (note the author had already identified and pointed out the mistake by Maskin and Dasgupta, and was working on the label-independent notions at the time of coming across Noah), and who answered a number of questions from the author, many of which were rather stupid in hindsight. Note the author purposefully omitted Noah from the acknowledgements in previous versions of this paper uploaded to the arXiv due to the author's uncertainty about whether the author would be considered a crank and dismissed without justification;
	    		\item Asaf Plan who identified a mistake with the proof of Theorem \ref{basicsymequivthm}, which the author has since corrected (the author apologises that he forgot to acknowledge Asaf in at least one version uploaded to the arXiv since Asaf pointed the mistake out);
	    		\item Jeremy Sumner for proof-reading early drafts with a number of useful corrections;
	    		\item Frigyes Podmaniczky who woke the author up to the annoying issues with the historical convention of listing players 1 and 2 along the rows and columns for games with $n \geq 3$ players;
	    		\item several anonymous referees whose useful comments have helped improve the exposition of this paper along with suggesting a number of relevant references; and
	    		\item Martin Osborne, whose {\tt sgamevar} {\LaTeX} style file, see \cite{SGAMESTY}, was used for the example games throughout this paper.
    		\end{enumerate} 
    
    \nocite{GAP4} 
    \bibliographystyle{plain}
    \bibliography{references}

\begin{thebibliography}{10}

\bibitem{ahern2018institutional}
Kathy Ahern.
\newblock Institutional betrayal and gaslighting.
\newblock {\em The Journal of perinatal \& neonatal nursing}, 32(1):59--65,
  2018.

\bibitem{aronowitz2009human}
Alexis~A Aronowitz.
\newblock {\em Human trafficking, human misery: The global trade in human
  beings}.
\newblock Greenwood Publishing Group, 2009.

\bibitem{arrow1950difficulty}
Kenneth~J. Arrow.
\newblock A difficulty in the concept of social welfare.
\newblock {\em Journal of political economy}, 58(4):328--346, 1950.

\bibitem{arrow2012social}
Kenneth~J Arrow.
\newblock {\em Social choice and individual values}, volume~12.
\newblock Yale university press, 2012.

\bibitem{brandt2009symmetries}
Felix Brandt, Felix Fischer, and Markus Holzer.
\newblock Symmetries and the complexity of pure {N}ash equilibrium.
\newblock {\em Journal of {C}omputer and {S}ystem {S}ciences}, 75(3):163--177,
  2009.

\bibitem{cao2018symmetric}
Zhigang Cao and Xiaoguang Yang.
\newblock Symmetric games revisited.
\newblock {\em Mathematical Social Sciences}, 95:9--18, 2018.

\bibitem{cao2019ordinally}
Zhigang Cao and Xiaoguang Yang.
\newblock Ordinally symmetric games.
\newblock {\em Operations Research Letters}, 47(2):127--129, 2019.

\bibitem{chantler2009forced}
Khatidja Chantler, Geetanjali Gangoli, and Marianne Hester.
\newblock Forced marriage in the {UK}: Religious, cultural, economic or state
  violence?
\newblock {\em Critical social policy}, 29(4):587--612, 2009.

\bibitem{CRVWSym}
Shih-Fen Cheng, Daniel~M Reeves, Yevgeniy Vorobeychik, and Michael~P Wellman.
\newblock Notes on equilibria in symmetric games.
\newblock In {\em Proceedings of the 6th International Workshop On Game
  Theoretic And Decision Theoretic Agents}, pages 71--78. GTDT, Research
  Collection School of Information Systems, 2004.

\bibitem{curato2019democracy}
Nicole Curato.
\newblock {\em Democracy in a time of misery: From spectacular tragedies to
  deliberative action}.
\newblock Oxford University Press, 2019.

\bibitem{DMaskin}
Partha Dasgupta and Eric Maskin.
\newblock The existence of equilibrium in discontinuous economic games, {I}:
  Theory.
\newblock {\em The Review of Economic Studies}, 53(1):18--25, 1986.

\bibitem{delli2002internet}
Michael~X Delli~Carpini and Scott Keeter.
\newblock The internet and an informed citizenry.
\newblock {\em Departmental Papers (ASC)}, page~2, 2002.

\bibitem{dungey2005contagion}
Mardi Dungey and Demosthenes Tambakis.
\newblock {\em Identifying {I}nternational {F}inancial {C}ontagion: {P}rogress
  and {C}hallenges}.
\newblock Oxford University Press, 2005.

\bibitem{EastHam}
James East and Nicholas Ham.
\newblock Lattice paths and submonoids of $\mathbb{Z}^2$.
\newblock {\em Preprint}, 2018, \href{http://arxiv.org/abs/1811.05735}{\tt
  arXiv:1811.05735}.

\bibitem{farago1941german}
Ladislas Farago.
\newblock German psychological warfare.
\newblock 1941.

\bibitem{IsoComplexity}
Joaquim Gabarr{\'o}, Alina Garc{\'\i}a, and Maria Serna.
\newblock The complexity of game isomorphism.
\newblock {\em Theoretical Computer Science}, 412(48):6675--6695, 2011.

\bibitem{gelbaum1959symmetric}
Bernard~R. Gelbaum.
\newblock Symmetric zero-sum n-person games.
\newblock {\em Annals of Mathematics Studies}, 40:95--109, 1959.

\bibitem{GoforthRobinson}
David Goforth and David Robinson.
\newblock {\em Topology of $2 \times 2$ Games}.
\newblock Routledge, 2005.

\bibitem{gold2011unit}
Hal Gold.
\newblock {\em Unit 731: Testimony}.
\newblock Tuttle Publishing, 2011.

\bibitem{GAP4}
The~GAP Group.
\newblock {GAP} -- {G}roups, {A}lgorithms, and {P}rogramming, {V}ersion 4.7.0,
  2013, {URL} \newline\href{https://www.gap-system.org}{\tt gap-system.org}.

\bibitem{ham2011honoursthesis}
Nicholas Ham.
\newblock Some {S}tructural {P}roperties of {F}inite {N}ormal {F}orm {G}ames.
\newblock {\em Honours thesis, University of Tasmania}, 2011.

\bibitem{ham2018arxivversion}
Nicholas Ham.
\newblock {N}otions of {S}ymmetry for {F}inite {S}trategic-{F}orm {G}ames.
\newblock {\em Preprint}, 2011,
  \newline\href{http://arxiv.org/abs/1311.4766}{\tt arXiv:1311.4766}.

\bibitem{GLCode}
Nicholas Ham.
\newblock {\em Code for generating finite strategic-form games}, December 2018,
  {URL}
  \newline\href{https://gitlab.com/n-ham-paper-files/generate-strategic-form-game}{\tt
  https://gitlab.com/n-ham-paper-files/generate-strategic-form-game}.

\bibitem{harriss2006poverty}
Barbara Harriss-White.
\newblock Poverty and capitalism.
\newblock {\em Economic and Political Weekly}, pages 1241--1246, 2006.

\bibitem{HarsanyiSelten}
John Harsanyi and Reinhard Selten.
\newblock {\em A General Theory of Equilibrium Selection in Games}.
\newblock MIT Press, 1988.

\bibitem{hefti2017equilibria}
Andreas Hefti.
\newblock Equilibria in symmetric games: Theory and applications.
\newblock {\em Theoretical Economics}, 12(3):979--1002, 2017.

\bibitem{hesse2011young}
Morten Hesse and S{\'e}bastien Tutenges.
\newblock Young tourists visiting strip clubs and paying for sex.
\newblock {\em Tourism Management}, 32(4):869--874, 2011.

\bibitem{hofbauer2002differential}
Josef Hofbauer and Gerhard Sorger.
\newblock A differential game approach to evolutionary equilibrium selection.
\newblock {\em International Game Theory Review}, 4(01):17--31, 2002.

\bibitem{hulpke2005constructing}
Alexander Hulpke.
\newblock Constructing transitive permutation groups.
\newblock {\em Journal of {S}ymbolic {C}omputation}, 39(1):1--30, 2005.

\bibitem{hulpketransgrp}
Alexander Hulpke, John Cannon, and Derek Holt.
\newblock Transgrp - {L}ibrary of {T}ransitive {G}roups.
\newblock {\em GAP manual}, year unknown.

\bibitem{hulpke2016connected}
Alexander Hulpke, David Stanovsk{\`y}, and Petr Vojt{\v{e}}chovsk{\`y}.
\newblock Connected quandles and transitive groups.
\newblock {\em Journal of Pure and Applied Algebra}, 220(2):735--758, 2016.

\bibitem{laffont2009theory}
Jean-Jacques Laffont and David Martimort.
\newblock {\em The theory of incentives: the principal-agent model}.
\newblock Princeton university press, 2009.

\bibitem{laffont1993theory}
Jean-Jacques Laffont and Jean Tirole.
\newblock {\em A theory of incentives in procurement and regulation}.
\newblock MIT press, 1993.

\bibitem{linebarger2015psychological}
Paul~MA Linebarger.
\newblock {\em Psychological warfare}.
\newblock Pickle Partners Publishing, 2015.

\bibitem{margolis2003misery}
Joshua~D Margolis and James~P Walsh.
\newblock Misery loves companies: Rethinking social initiatives by business.
\newblock {\em Administrative science quarterly}, 48(2):268--305, 2003.

\bibitem{Mas1995microeconomic}
Andreu Mas-Colell, Michael~Dennis Whinston, Jerry~R Green, et~al.
\newblock {\em Microeconomic theory}, volume~1.
\newblock Oxford university press New York, Pages 9, 46-48, 1995.

\bibitem{minkowitz2006united}
Tina Minkowitz.
\newblock The united nations convention on the rights of persons with
  disabilities and the right to be free from nonconsensual psychiatric
  interventions.
\newblock {\em Syracuse J. Int'l L. \& Com.}, 34:405, 2006.

\bibitem{NashNCG}
John Nash.
\newblock Non-cooperative games.
\newblock {\em The Annals of Mathematics, Second Series}, 54(2):286--295, 1951.

\bibitem{SGAMESTY}
Martin~J. Osborne.
\newblock {\LaTeX} style files {\tt sgame.sty} and {\tt sgamevar.sty} for
  drawing strategic-form games, {V}ersion 2.15, 2012, {URL}
  \url{https://www.economics.utoronto.ca/osborne/latex/}.

\bibitem{ouattara1998forced}
Mariam Ouattara, Purna Sen, and Marilyn Thomson.
\newblock Forced marriage, forced sex: the perils of childhood for girls.
\newblock {\em Gender \& Development}, 6(3):27--33, 1998.

\bibitem{peleg1999canonical}
Bezalel Peleg, Joachim Rosenm{\"u}ller, and Peter Sudh{\"o}lter.
\newblock The canonical extensive form of a game form: Symmetries.
\newblock In {\em Current Trends in Economics. Studies in Economic Theory},
  volume~8, pages 367--387. Springer, Berlin, Heidelberg, 1999.

\bibitem{shapley1960symmetric}
Lloyd~S. Shapley.
\newblock Symmetric games.
\newblock {\em RAND Corporation}, 1960.

\bibitem{shelley2010human}
Louise Shelley.
\newblock {\em Human trafficking: A global perspective}.
\newblock Cambridge University Press, 2010.

\bibitem{NoahXE}
Noah Stein.
\newblock Exchangeable {E}quilibria.
\newblock {\em Doctoral Thesis, MIT}, 2011.

\bibitem{Steinhauscakecutting}
Hugo Steinhaus.
\newblock The {P}roblem of {F}air {D}ivision. {R}eport of the {W}ashington
  {M}eeting, {S}eptember 6-18, 1947.
\newblock {\em Econometrica}, 16(1):101--104, 1948.

\bibitem{sudholter2000canonical}
Peter Sudh{\"o}lter, Joachim Rosenm{\"u}ller, and Bezalel Peleg.
\newblock The canonical extensive form of a game form: Part {II}.
  {R}epresentation.
\newblock {\em Journal of Mathematical Economics}, 33(3):299--338, 2000.

\bibitem{sweet2019sociology}
Paige~L Sweet.
\newblock The sociology of gaslighting.
\newblock {\em American Sociological Review}, 84(5):851--875, 2019.

\bibitem{tohme2019structural}
Fernando~A Tohm{\'e} and Ignacio~D Viglizzo.
\newblock Structural relations of symmetry among players in strategic games.
\newblock {\em International Journal of General Systems}, 48(4):443--461, 2019,
  \href{http://arxiv.org/abs/1712.04563}{\tt arXiv:1712.04563}.

\bibitem{tucker1962comb}
Albert~W. Tucker.
\newblock Combinatorial equivalence of fair games.
\newblock In {\em Recent Advances in Game Theory: Papers Delivered at a Meeting
  of the Princeton University Conference, Ocotber [sic] 4-6, 1961}, volume~29,
  page 277, 1962.

\bibitem{vester2012symmetric}
Steen Vester.
\newblock Symmetric {N}ash equilibria.
\newblock {\em Masters thesis, Ecole Normale Superieure de Cachan}, 2012.

\bibitem{VNM}
John von Neumann and Oskar Morgenstern.
\newblock {\em Theory of Games and Economic Behaviour}.
\newblock Princeton University Press, New Jersey, 1944.

\bibitem{walker1966critique}
Jack~L Walker.
\newblock A critique of the elitist theory of democracy.
\newblock {\em The American Political Science Review}, 60(2):285--295, 1966.

\bibitem{WikiSGRV}
Wikipedia.
\newblock Symmetric games: {R}evision history.
\newblock Accessed: 2018-12-06, {URL}
  \newline\url{https://en.wikipedia.org/w/index.php?title=Symmetric_game&action=history}.

\bibitem{williams1989unit}
Peter Williams and David Wallace.
\newblock {\em Unit 731: Japan's secret biological warfare in World War II}.
\newblock Free Press New York, 1989.

\end{thebibliography}
\end{document}